\documentclass[journal]{IEEEtran}      % Use this line for a4
\IEEEoverridecommandlockouts                              % This command is only
\usepackage{mathrsfs}

\usepackage{subfigure}
\usepackage{graphicx}
\usepackage{amsmath,amsfonts,amssymb}
\usepackage{array}

\def \R {\mathbb{R}}
\def \Y {\mathcal{Y}}
\def \N {\mathbb{N}}

\def \bv {\bar{v}}

\newtheorem{definition}{Definition}

\newtheorem{proposition}{Proposition}

\newtheorem{theorem}{Theorem}
\renewcommand{\thetheorem}{\arabic{section}.\arabic{theorem}}
\newtheorem{lemma}{Lemma}

\newtheorem{corollary}{Corollary}

\newtheorem{remark}{Remark}

\title{\LARGE \bf
Multivariate Feedback Particle Filter via F-divergence and the Well-posedness of Its Admissible Control Input}

\author{%~\IEEEmembership{Member,~IEEE,}
        Xue~Luo,~\IEEEmembership{senior member,~IEEE},
     \thanks{This work is financially supported by National Natural Science Foundation of China (NSFC, grant no. 11871003, 11501023, 11471184) and the Fundamental
Research Funds for the Central Universities (grant no.
YWF-18-BJ-J-238).}% <-this % stops a space
\thanks{X. Luo is with School of Mathematics and Systems Science, Beihang University, Beijing, P. R. China        {\tt\small xluo@buaa.edu.cn}}
}

\def \Y {\mathcal{Y}}
\def \N {\mathbb{N}}
\def \R {\mathbb{R}}
\def \Rd {\mathbb{R}^d}

\def \L {\mathcal{L}}
\def \E {\mathbb{E}}

\begin{document}

\maketitle
\thispagestyle{empty}

\begin{abstract}
	In this paper, we shall first derive the admissible control input of the multivariate feedback particle filter (FPF) by minimizing the f-divergence of the posterior conditional density function and the empirical conditional density of the controlled particles. On the contrast, in the original derivation \cite{YMM}, a special f-divergence, Kullback-Leibler (K-L) divergence, is used in the 1-dimensional nonlinear filtering problems. We show that the control input is invariant under the f-divergence class. That is, the control input satisfies exactly the same equations as those obtained by minimizing K-L divergence, no matter what f-divergence in use. In the latter half of this paper, we show the existence and uniqueness of the control input under suitable regular conditions. We confirm that the explicit expression of the control input given in \cite{YLMM} is the only admissible one.
\end{abstract}

\section{Introduction}

Nonlinear filtering (NLF) or called nonlinear estimation is to give the state/signal $X_t$ a ``proper" estimation based on the observation history $\{Z_s,0\leq s\leq t\}$ in some sense, say the conditional expectation $\E[X_t|\mathcal{Z}_t]$, $\mathcal{Z}_t:=\sigma(\{Z_s,\,0\leq s\leq t\})$. The most famous Kalman filter \cite{Ka}, Kalman-Bucy filter \cite{KB} yields the optimal estimation if the problem is linear and the initial density is Gaussian. Unfortunately, as far as we know, there is no such effective algorithms for NLF problems. There are lots of Kalman filter's derivatives which can obtain suboptimal estimations for NLF problems, but far from satisfactory.  We refer the approaches that only aim to obtain the approximation of interested statistical quantities, say expectation, variance  etc., as local approaches, while those consider to compute/approximate the posterior distribution of the hidden Markov process $X_t$, given the history of observation is called global approaches. The author of this paper wrote a survey on the approaches for NLF problems and carefully discussed their advantage/disadvantage, see \cite{Lu}. 

One global approach is to numerically solving the Kushner-Stratonovich equation  \cite{K} or Duncan-Mortensen-Zakai equation \cite{D,Mo,Z}. In this direction, there are wide range of literatures including the splitting-up method \cite{BGR}, $S^3$-algorithm \cite{LMR}, on- and off-line algorithm \cite{YY,LY1,LY2}, etc. Nevertheless, the computation load is extremely heavy if the state/signal is high-dimensional. Besides these partial differential equation based algorithm, the most popular global approach is the so-called particle filter (PF) \cite{AMG,CD,DFG}. The PF is a simulation-based algorithm, which approximates the posterior distribution by the empirical distribution of the particle population $\{X_t^i\}_{i=1,\cdots,N}$. A common remedy to avoid particle impoverish and degeneracy in the traditional PF is to vigor the particles by resampling according to the importance weight at each time step. After the proper reampling strategy, the PF can propagate the posterior distribution with accuracy improving by increasing the number of the particles \cite{CD}. Nevertheless, the choice of the importance weight is crucial, problem-dependent and with no universal guidelines.

Recently, Yang et. al. \cite{YMM} introduced a control-oriented PF, called feedback particle filter (FPF), for the scalar case, i.e. the dimensional of the state/signal and observation process both are one. Later, \cite{YLMM} extends this algorithm to multivariate case without detailed derivation. The FPF is motivated by mean-field optimal control techniques \cite{LL,HCM}. Let us consider the NLF problem in the form:
\begin{align}\label{eqn-NLF continuous}
	dX_t=&a(X_t)dt+\sigma_BdB_t\\\label{eqn-observation continuous}
	dZ_t=&h(X_t)dt+dW_t,
\end{align}
where $X_t\in\Rd$ is the state at time $t$, $Z_t\in\R^m$ is the observation process, $a(\cdot)$, $h(\cdot)$ are functions of $X_t$ and $\{B_t\}$, $\{W_t\}$ are mutually independent Wiener processes of appropriate dimensions. In FPF, they model the $i$th particle evolves according the controlled system
\begin{align}\label{eqn-particle controlled}
	dX_t^i=a(X_t^i)dt+\sigma_BdB_t^i+dU_t^i,
\end{align}
where $dU_t^i$ is the control input of the $i$th particle, and $\{B_t^i\}$ are also mutually independent standard Wiener process. The aim of the FPF is to choose the appropriate control input for each particle such that the empirical distribution approximates the conditional distribution of $X_t$ for large number of particles. In \cite{YMM,YLMM}, the optimal control input is obtained by minimizing the Kullback-Leibler (K-L) divergence as the cost function. 

In this paper, we shall discuss two natural questions related to the multivariate FPF. On the one hand, in probability theory the K-L divergence is only one member in the category called f-divergence, which measures the difference between two probability distributions. These f-divergences were introduced and studied independently by Csisz\'ar \cite{C}, Morimoto \cite{M} and Ali, et. al. \cite{AS} and are sometimes known as Csisz\'ar f-divergences, Csisz\'ar-Morimoto divergences or Ali-Silvey distances. Thus, a natural question is raised: if the other f-divergences are used as the cost function in FPF, shall we obtain different control input from that obtained by using the K-L divergence? We answered this question in the scalar case, i.e. $d=m=1$, in \cite{LW} that the control input is independent of the choice of f-divergence. When it comes to the multivariate case, it is not trivial. As one will see in this paper, the trivial identity \eqref{eqn-nabla=0} for $d=1$, which has to be shown rigorously for $d\geq2$, see Proposition \ref{lem-for any d} in section III.A. 

On the other hand, the derivation of the multivariate FPF in \cite{YLMM} is too informal for the readers to suspect that the control input is just an ananalogue of the one in the scalar case, without any detailed derivation, let alone the discussion of the existence and uniqueness of the control input in what sense. In section III.B, we patch the detailed derivation for the equations which the control input should satisfy for the case $d\geq1$, $m=1$. Consequently, the consistency can be shown rigorously based on the equations derived, rather than on the analoguous control input ``guessed" in \cite{YLMM}. For the most general case $m\geq1$, we point out that our derivation should also work but with more involved computations and notations. In section IV, we established the well-posedness of the control input in appropriate function space under certain conditions. Thus, the control input given in \cite{YLMM} has been checked to be admissible (Definition \ref{def-admissible}), so as to be unique in the suitable function space. The conclusions are arriven in the end.

\section{Preliminaries}

The precise formulation begins with continuous time model \eqref{eqn-NLF continuous} with sampled observations:
\begin{equation}
	Y_{t_n}=h(X_{t_n})+W^{\triangle}_{t_n},
\end{equation}
where $\Delta t:=t_{n+1}-t_n$ and $\{W^{\Delta t}_{t_n}\}$ is i.i.d. and drawn from $\mathcal{N}(0,\frac1{\Delta t})$, the Gaussian with zero mean, $\frac1{\Delta t}$ variance. The observation history is denoted as $\Y_n:=\sigma\{Y_{t_k}:\ k\leq n, k\in\N\}$. We follow the same notations as in \cite{YMM, YLMM}. Let us denote the conditional distributions:
\begin{enumerate}
	\item $p_n^*$ and $p_n^{*-}$: the conditional distribution of $X_{t_n}$ given $\Y_n$ and $\Y_{n-1}$, respectively.
	\item $p_n$ and $p_n^-$: the conditional distribution of the $i$th particle $X_{t_n}^i$ given $\Y_n$ and $\Y_{n-1}$, respectively.
\end{enumerate}
These distributions evolve according to the recursion
\begin{align*}
	p_n^*=\mathcal{P}^*(p_{n-1}^*,Y_{t_n}),\quad 
	p_n=\mathcal{P}(p_{n-1},Y_{t_n}).
\end{align*}
The mappings $\mathcal{P}^*$ and $\mathcal{P}$ can be  decomposed into two parts. The first part is identical for each of these mappings: the transformation that takes $p_{n-1}$ to $p^-_n$ coincides with the mapping from $p^*_{n-1}$ to $p^{*-}_n$.In each case it is defined by the Kolmogorov forward equation (KFE) associated with the diffusion on $[t_{n-1},t_n)$.
 
The second part of the mapping is the updating that takes $p^{*-}_n$ to $p^*_n$ by synchronizing the observation data $Y_{t_n}$, which is obtained according to the Bayes' rule. Given the observation $Y_{t_n}$ made at time $t=t_n$
\begin{equation}\label{eqn-1}
p^*_n(s)=\frac{p^{*-}_n(s)\cdot p_{Y|X}(Y_{t_n}|s)}{p_Y(Y_{t_n})},\quad s\in \Rd,
\end{equation}
where $p_Y$ denotes the probability density function (pdf) for $Y_{t_n}$, and $ p_{Y|X}(Y_{t_n} |s)$ denotes the conditional distribution of $Y_{t_n}$ given $X_{t_n}=s$. In the case that the observation noise is Gaussian, we have
\begin{equation}\label{eqn-2}
p_{Y|X}(Y_{t_n}|s)=\frac{1}{\sqrt{\frac{2\pi}{\Delta t}}} \textup{exp}\left(-\frac{(Y_{t_n}-h(s))^2}{\frac{2}{\Delta t}}\right)
\end{equation}
The operator $\mathcal{P}^*$ is the composition of KFE and \eqref{eqn-2}.

The updating from $p^-_n$ to $p_n$ is not due to the Bayes' rule, but depends on the control input $dU^i_{t_n}$ in \eqref{eqn-particle controlled}. In the discrete setting, at time $t=t_n$, we seek the control input $v_n^i=K(X_{t_n}^i,t_n)\Delta z+u(X_{t_n}^i,t_n)\Delta t$, which is the discrete counterpart of $dU_{t_n}^i=K(X_t^i,t)dZ_t+u(X_t^i,t)dt$ at time $t=t_n$. We shall restrict ourselves to find the control input in the {\em admissible} class. We will suppress the superscript $i$ and the subscript $n$ in $v_n^i$, and write $K=K(X_{t_n}^i,t_n)$, $u=u(X_{t_n}^i,t_n)$ for short, if there is no confusion.

\begin{definition} [Admissible Input]\label{def-admissible}
 The control sequence $v_n=K\Delta z+u\Delta t$ is called admissible, if for each $n$
	\begin{enumerate}
	\item $K\in H^l(\R^d;p)$, $l\geq\lfloor\frac d2\rfloor+1$ and $u\in L^2(\R^d;p)$, where $H^l(\R^d;p)$ is the weighted Sobolev space with its norm defined as
\begin{align*}
	||\circ||_{H^l(\R^d;p)}^2=\sum_{i=0}^l||\nabla^i\circ||_{L^2(\R^d;p)}^2,
\end{align*} 
with 
\begin{equation}\label{eqn-nablasigmak}
	\nabla_{\sigma_i}^i\circ:=\frac{\partial^i\circ}{\partial x_{\sigma_i(1)}\cdots\partial x_{\sigma_i(i)}},
\end{equation}
and the convention that $\nabla^0=Id$, the identity mapping, where the norm is defined as
\begin{align}\label{eqn-norm of L2}\notag
	&||\nabla^i\circ||^2_{L^2(\R^d;p)}\\
	=&\sum_{j=1}^d\sum_{\sigma_i\in\{1,\cdots,d\}^l}\int_{\R^d}|\nabla^i_{\sigma_i}\circ_j|^2pdx,
\end{align}
for $\circ=(\circ_1,\cdots,\circ_d)\in\R^d$.
	\item $I+(\nabla v^T_n)^T$ is invertible for all $x$, where $I$ is the identity matrix, $\circ^T$ is the transpose of the matrix $\circ$ and $\nabla v^T=
\begin{pmatrix}
	\frac\partial{\partial x_1}\\
	\vdots\\
	\frac\partial{\partial x_d}
\end{pmatrix}(v_1,\cdots,v_d)=
\begin{pmatrix}
	\frac{\partial v_1}{\partial x_1}\ \cdots\ \frac{\partial v_d}{\partial x_1}\\
	\vdots\ \ddots\ \vdots\\
	\frac{\partial v_1}{\partial x_d}\ \cdots\ \frac{\partial v_d}{\partial x_d}
\end{pmatrix}$, which is the transpose of the Jacobian matrix $\frac{\partial v}{\partial x}$.
\end{enumerate}

 Under the assumption that $I+(\nabla v^T)^T$ is invertible for all $x$, the updating from $p_n^-$ to $p_n$ is 
\begin{equation}\label{eqn-3}
p_n(s)=\frac{p^-_n(x)}{|I+\nabla v^T(x)|}, 
\end{equation}
where $s=x+v(x)$, and $|\circ|$ represents the determinant of the matrix $\circ$.
\end{definition}

Let us denote $\hat{p}^*=\mathcal{P}^*(p_{n-1},Y_{t_n})$. We should choose the control input $v_n^i$ such that the difference between $\hat{p}^*_n$ and $p_n$ as small as possible at every time step $t_n$. 

In 1960s, the evolution equation for $p^*$ has been derived by Kushner \cite{K}:
\begin{equation}\label{eqn-Kushner}
	dp^*=\L^*p^*dt+(h-\hat{h}_t)(dZ_t-\hat{h}_tdt)p^*,
\end{equation}
where $\hat{h}_t=\int hp^*dx$ and 
\begin{equation}\label{eqn-L^*}
	\L^*p^*:=-\nabla^T(p^*a)+\frac12\sum_{i,j=1}^d\frac{\partial^2}{\partial x_i\partial x_j}(p^*[\sigma_B\sigma_B^T]_{ij}).
\end{equation} 
This is the so-called the Kushner-Stratonovich equation. Moreover, the propagation of the particle's conditional distribution is described by the Kolmogorov forward equation (KFE) \cite{O}:
\begin{align}\label{eqn-KFE}\notag
	dp=&\L^*pdt-\nabla^T(pK)dZ_t-\nabla^T(pu)dt\\
	&+\frac12\sum_{i,j=1}^d\frac{\partial^2}{\partial x_i\partial x_j}(p[KK^T]_{ij})dt,
\end{align}
where $\L^*$ is defined in \eqref{eqn-L^*}. The derivation of this KFE in scalar case can be found in Proposition 3.1, \cite{YMM}. The multivariate case is straight-forward.

\section{Multivariate feedback particle filter}

\subsection{Invariance of control input via F-divergence}

In this section, we measure the difference from distribution $p_1$ to $p_2$ by $f$-divergence defined as
\begin{align*}
	D_f(p_1||p_2)=\int_{\Rd}p_2(s)f\left(\frac{p_1(s)}{p_2(s)}\right)ds.
\end{align*}
With different choice of $f$, f-divergence can become Kullback-Leibler (K-L) divergence, total variation distance, Hellinger distance, etc. The K-L divergence is a special case of $f$-divergence with $f(s)=s\log s$. In \cite{YMM,YLMM}, the control input is obtained by minimizing the K-L divergence from $p_n$ to $\hat{p}_n^*$.

In this section, we shall determine the control input $v$ by minimizing the f-divergence from $\hat{p}^*_n$ and $p_n$. Although it is well-known that the $f$-divergence is not symmetry, i.e. $D_f(p_n||\hat{p}_n^*)\neq D_f(\hat{p}_n^*||p_n)$, the control input $v$ obtained by minimizing either one is exactly the same. Nevertheless, the computation of $D_f(\hat{p}_n^*||p_n)$ is much easier. Thus, we use $D_f(\hat{p}_n^*||p_n)$ in the derivation:
\begin{align}\label{eqn-Df}\notag
	&D_f(\hat{p}_n^*||p_n)\\
	\overset{\eqref{eqn-1},\eqref{eqn-3}}	=&\int_{\Rd} \frac{p_n^-(x)}{|I+\nabla v^T|}f\left(\frac{\frac{p^-_n(s)p_{Y|X}(Y_{t_n}|s)}{p_Y(Y_{t_n})}}{\frac{p_n^-(x)}{|I+\nabla v^T|}}\right)ds\\\notag
	=&\int_{\Rd} p_n^-(x)f\left(\frac{p_n^-(s)p_{Y|X}(Y_{t_n}|s)|I+\nabla v^T|}{p_n^-(x)p_Y(Y_{t_n})}\right)dx,
\end{align}
where $s=x+v(x)$. Let us denote the integrand in \eqref{eqn-Df} as
\begin{align}\label{eqn-Lxvnablav}\notag
	&\L(x,v,\nabla v^T)\\
	:=&p_n^-(x)f\left(\frac{p_n^-(s)p_{Y|X}(Y_{t_n}|s)|I+\nabla v^T|}{p_n^-(x)p_Y(Y_{t_n})}\right).
\end{align}
For simplicity of notation, let us denote the argument of $f$ in \eqref{eqn-Lxvnablav} as  
\begin{equation}\label{eqn-xi}
	\xi=\frac{p_n^-(x+v)p_{Y|X}(Y_{t_n}|x+v)|I+\nabla v^T|}{p_n^-(x)p_Y(Y_{t_n})}.
\end{equation}

It is well-known that the minimizer of the functional $D_f(\hat{p}_n^*||p_n)$ can be obtained by solving the corresponding Euler-Lagrange (E-L) equation:
\begin{align}\label{eqn-EL}
	\left(\frac{\partial\L}{\partial v}\right)^T=\nabla^T_x\left(\frac{\partial \L}{\partial (\nabla v^T)}\right),
\end{align}
where taking the derivative with respect to a vector or a matrix has been defined properly in matrix calculus, see Chapter 9, \cite{L}. The notation $\nabla_x^T(A)$ in \eqref{eqn-EL} is $\nabla_x^T(A)=\begin{pmatrix}
	\sum_{j=1}^d\frac{\partial A_{1j}}{\partial x_j},\cdots,\sum_{j=1}^d\frac{\partial A_{dj}}{\partial x_j}
\end{pmatrix}$, if $A$ is a $d\times d$ matrix and $A_{ij}$ is the $ij$th entry of $A$. 
 The left-hand side of \eqref{eqn-EL} equals:
\begin{align}\label{eqn-lhs of EL}\notag
	&\left(\frac{\partial\L}{\partial v}\right)^T\\\notag
	=&\frac{f'(\xi)|I+\nabla v^T|}{p_Y(Y_{t_n})}\,\nabla^T\left(p_n^-(x+v)p_{Y|X}(Y_{t_n}|x+v)\right)\\
	=&\frac{f'(\xi)|I+\nabla v^T|}{p_Y(Y_{t_n})}\\\notag
	&\cdot\nabla_x^T\left(p_n^-(x+v)p_{Y|X}(Y_{t_n}|x+v)\right)(I+\nabla v^T)^{-T},
\end{align}
where the notation $A^{-T}=(A^{-1})^T$, while its right-hand side is
\begin{align}\label{eqn-rhs of EL}\notag
	&\nabla^T_x\left(\frac{\partial \L}{\partial (\nabla v^T)}\right)\\\notag
	=&\nabla_x^T\left[f'(\xi)\frac{p_n^-(x+v)p_{Y|X}(Y_{t_n}|x+v)}{p_Y(Y_{t_n})}\frac{\partial |I+\nabla v^T|}{\partial(\nabla v^T)}\right]\\\notag
	=&\nabla_x^T\left[f'(\xi)\frac{p_n^-(x+v)p_{Y|X}(Y_{t_n}|x+v)}{p_Y(Y_{t_n})}\right.\\\notag
	&\phantom{\nabla_x^T[]}\cdot\left.|I+\nabla v^T|(I+\nabla v^T)^{-T}\right]\\\notag
	=&p_n^-(x+v)p_{Y|X}(Y_{t_n}|x+v)\\\notag
	&\phantom{+}\cdot\nabla_x^T\left[f'(\xi)\frac{|I+\nabla v^T|}{p_Y(Y_{t_n})}\right](I+\nabla v^T)^{-T}\\\notag
	&+f'(\xi)\frac{|I+\nabla v^T|}{p_Y(Y_{t_n})}\nabla_x^T[p_n^-(x+v)p_{Y|X}(Y_{t_n}|x+v)]\\\notag
	&\phantom{+}\cdot(I+\nabla v^T)^{-T}\\\notag
	&+p_n^-(x+v)p_{Y|X}(Y_{t_n}|x+v)f'(\xi)\frac{|I+\nabla v^T|}{p_Y(Y_{t_n})}\\
	&\phantom{+}\cdot\nabla_x^T\left[(I+\nabla v^T)^{-T}\right],
\end{align}
where the second equality in \eqref{eqn-rhs of EL} follows from the fact that $\frac{\partial |A|}{\partial A}=|A|A^{-T}$. Notice that the second term on the right-hand side of \eqref{eqn-rhs of EL} cancels out with \eqref{eqn-lhs of EL}. Therefore, we obtain from \eqref{eqn-EL} that
\begin{align}\label{eqn-EL-f}\notag
	0=&\nabla_x^T\left(f'(\xi)|I+\nabla v^T|\right)(I+\nabla v^T)^{-T}\\\notag
	&+f'(\xi)|I+\nabla v^T|\nabla_x^T\left[(I+\nabla v^T)^{-T}\right]\\\notag
	=&\nabla_x^T[f'(\xi)]|I+\nabla v^T|(I+\nabla v^T)^{-T}\\
	&+f'(\xi)\nabla_x^T\left[|I+\nabla v^T|(I+\nabla v^T)^{-T}\right],
\end{align}
after dividing by $\frac{p_n^-(x+v)p_{Y|X}(Y_{t_n}|x+v)}{p_Y(Y_{t_n})}$ throughout.

We can show that actually the second term on the right-hand side of \eqref{eqn-EL-f} vanishes. 
\begin{proposition}\label{lem-for any d}
	For $x\in\R^d$, $d\geq1$, the identity
	\begin{equation}\label{eqn-nabla=0}
		\nabla_x^T\left[|I+\nabla v^T|(I+\nabla v^T)^{-T}\right]\equiv0
	\end{equation}
holds.
\end{proposition}

The detailed proof has been appended in Appendix A for interested readers.

Hence, the control input $v$ satisfies
\begin{equation}\label{eqn-EL-f1}
	0=\nabla_x^T[f'(\xi)]|I+\nabla v^T|(I+\nabla v^T)^{-T},
\end{equation}
where $\xi$ is defined in \eqref{eqn-xi}. Let us take a look at the term $\nabla_x^Tf'(\xi)$:
\begin{align}\label{eqn-nabla f'}\notag
	&\nabla_x^Tf'(\xi)\\\notag
	=&f''(\xi)\nabla_x^T\left[\frac{p_n^-(x+v)p_{Y|X}(Y_{t_n}|x+v)|I+\nabla v^T|}{p_n^-(x)p_Y(Y_{t_n})}\right]\\
	=&\frac{f''(\xi)}{p_Y(Y_{t_n})[p_n^-(x)]^2}\\\notag
	&\cdot\left\{\nabla_x^T[p_n^-(x+v)p_{Y|X}(Y_{t_n}|x+v)]|I+\nabla v^T|p_n^-(x)\right.\\\notag
	&\phantom{\cdot[]}+p_n^-(x+v)p_{Y|X}(Y_{t_n}|x+v)|I+\nabla v^T|\\\notag
	&\phantom{\cdot[]+}\cdot\left[\textup{tr}\left((I+\nabla v^T)^{-1}\frac{\partial(\nabla v^T)}{\partial x_i}\right)\right]_{i=1,\cdots,d}p_n^-(x)\\\notag
	&\phantom{\cdot[]}\left.-p_n^-(x+v)p_{Y|X}(Y_{t_n}|x+v)|I+\nabla v^T|\nabla_x^Tp_n^-(x)\right\},
\end{align}
where $\textup{tr}(\circ)$ means the trace of the matrix $\circ$ and the second term on the right-hand side follows from the fact that $\frac{d|A(x)|}{dx}=|A|\textup{tr}\left(A^{-1}\frac{dA}{dx}\right)$, with $A$ being a matrix-valued function of $x$. 

If the terms in the brace on the right-hand side of \eqref{eqn-nabla f'} equals zero, then the control input is independent of what $f$-divergence we are using. We summarize this invariance in the following theorem:
\begin{theorem} 
The control input $v$ obtained by minimizing the f-divergence from $\hat{p}_n^*$ to $p_n$ is independent of the choice of f-divergence in use. No matter what $f$ is, the control input $v$ always satisfies the following equation:
\begin{align}\label{eqn-nabla f'=0}\notag
	0=&\nabla_x^T[p_n^-(x+v)p_{Y|X}(Y_{t_n}|x+v)]p_n^-(x)\\\notag
	&+p_n^-(x+v)p_{Y|X}(Y_{t_n}|x+v)\\\notag
	&\phantom{+}\cdot\left[\textup{tr}\left((I+\nabla v^T)^{-1}\frac{\partial(\nabla v^T)}{\partial x_i}\right)\right]_{i=1,\cdots,d}p_n^-(x)\\
	&-p_n^-(x+v)p_{Y|X}(Y_{t_n}|x+v)\nabla_x^Tp_n^-(x).
\end{align}
\end{theorem}
\begin{remark}
	Equation \eqref{eqn-nabla f'=0} holds for any $p_{Y|X}$. That is, the observation noise can be other type beyond the Gaussian. For example, the author of this paper and her co-work investigate the FPF for observation noise with Laplace distribution \cite{LW}.  
\end{remark}

\subsection{Consistency}

Let us back to the situation that the observation noise is Gaussian as in \cite{YLMM}, i.e. 
\begin{align*}
	&\nabla p_{Y|X}(Y_{t_n}|x+v)\\
	\overset{\eqref{eqn-2}}=&\nabla\left[\frac1{\sqrt{\frac{2\pi}{\triangle t}}}\exp\left(-\frac{(Y_{t_n}-h(x+v))^2}{\frac{2}{\triangle t}}\right)\right]\\
	=&p_{Y|X}(Y_{t_n}|x+v)(\triangle z-h(x+v)\triangle t)\nabla h(x+v),
\end{align*}
where $\triangle z=\frac{Y_{t_n}}{\triangle t}$. Thus, \eqref{eqn-nabla f'=0} becomes
\begin{align}\label{eqn-nabla f'=0-1}\notag
	0=&\nabla_x^Tp_n^-(x+v)p_n^-(x)\\\notag
	&+p_n^-(x+v)(\triangle z-h(x+v)\triangle t)\nabla_x^Th(x+v)p_n^-(x)\\\notag
	&+p_n^-(x+v)\\\notag
	&\phantom{+}\cdot\left[\textup{tr}\left((I+\nabla v^T)^{-1}\frac{\partial (\nabla v^T)}{\partial x_i}\right)\right]_{i=1,\cdots,d}p_n^-(x)\\\notag
	&-p_n^-(x+v)\nabla_x^Tp_n^-(x)\\
	&=:II_1+II_2+II_3+II_4,
\end{align}
after dividing by $p_{Y|X}(Y_{t_n}|x+v)$ throughout. For the conciseness of the notation, we shall suppress $p_n^-(x)$ as $p$ below, if no confusion will arise. We shall seek for the control input in the form $v=K\Delta z+u\Delta t$. The Taylor expansion around $x$ is applied to $II_1-II_4$ one-by-one:
\begin{align*}
	II_1
=&\nabla^Tp(x+v)p(I+\nabla v^T)^T\\
=&p\left[\nabla^Tp+v^T\nabla^2p+\frac12v^T[\nabla\otimes(\nabla^2p)](I\otimes v)\right]\\
	&\cdot(I+\nabla v^T)^T\\
=&p\left\{\nabla^Tp+K^T\nabla^2p\triangle z\right.\\
	&\phantom{p[]}\left.+\left[u^T\nabla^2p+\frac12K^T[\nabla^T\otimes(\nabla^2p)(I\otimes K)]\right]\triangle t\right\}\\
	&\phantom{p[]+}\cdot(I+\nabla K^T\triangle z+\nabla v^T\triangle t)^T\\
=&p\nabla^Tp+\left[pK^T\nabla^2p+p\nabla^Tp(\nabla K^T)^T\right]\triangle z\\
&+\left\{p\left[u^T\nabla^2p+\frac12K^T[\nabla^T\otimes(\nabla^2p)(I\otimes K)]\right]\right.\\
	&\phantom{+[]}\left.+p\nabla^Tp(\nabla u^T)^T+pK^T\nabla^2p(\nabla K^T)^T\right\}\triangle t,
\end{align*}
\begin{align*}
	II_2
=&p\left[p+\nabla^Tpv+\frac12v^T(\nabla^2p)v\right]\\
	&\cdot\left[\triangle z-(h+\nabla^Thv+\frac12v^T(\nabla^2h)v)\triangle t\right]\\
	&\cdot\left[\nabla^Th+v^T\nabla^2h+\frac12v^T[\nabla^T\otimes(\nabla^2h)](I\otimes v)\right]\\
	&\cdot\left[I+(\nabla K^T)^T\triangle z+(\nabla u^T)^T\triangle t\right]\\
=&p^2\nabla^Th\triangle z\\
	&+\left[-p^2h\nabla^Th+p\nabla^TpK\nabla^Th\right.\\
	&\phantom{+[]}\left.+p^2K^T\nabla^2h+p^2\nabla^Th(\nabla K^T)^T\right]\triangle t,
\end{align*}
\begin{align*}
	II_3
=&p\left[p+(\nabla^Tp)K\triangle z+\left((\nabla^Tp)u+\frac12K^T(\nabla^2p)K\right)\triangle t\right]\\
	&\cdot\left[\textup{tr}\left((I+\nabla v^T)^{-1}\frac\partial{\partial x_i}\left(\nabla K^T\triangle z+\nabla u^T\triangle t\right)\right)\right]_{i=1,\cdots,d}\\
	=&p\left[p+(\nabla^Tp)K\triangle z+\left((\nabla^Tp)u+\frac12K^T\nabla^2pK\right)\triangle t\right]\\
	&\cdot\left\{\textup{tr}\left[\frac{\partial(\nabla K^T)}{\partial x_i}\triangle z\right.\right.\\
	&\phantom{\cdot\textup{tr}aaaa}\left.\left.+\left(\frac{\partial(\nabla u^T)}{\partial x_i}-\nabla K^T\frac{\partial(\nabla K^T)}{\partial x_i}\right)\triangle t\right]\right\}_{i=1,\cdots,d}\\
	=&p^2\left[\textup{tr}\left(\frac{\partial}{\partial x_i}(\nabla K^T)\right)\right]_{i=1,\cdots,d}\triangle z\\
	&+\left\{p\nabla^TpK\left[\textup{tr}\left(\frac{\partial (\nabla K^T)}{\partial x_i}\right)\right]_{i=1,\cdots,d}\right.\\
	&\phantom{+aa}\left.+p^2\left[\textup{tr}\left(\frac{\partial(\nabla u^T)}{\partial x_i}-\nabla K^T\frac{\partial(\nabla K^T)}{\partial x_i}\right)\right]_{i=1,\cdots,d}\right\}\triangle t,
\end{align*}
and 
\begin{align*}
	II_4
=&-\left[p+\left(\nabla^Tp\right)K\triangle z\right.\\
	&\phantom{-aa}\left.+\left(\left(\nabla^Tp\right)u+\frac12K^T\left(\nabla^2p\right)K\right)\triangle t\right]\nabla^Tp\\
=&-p\nabla^Tp-\left(\nabla^Tp\right)K\nabla^Tp\triangle z\\
	&-\left(\left(\nabla^Tp\right)u+\frac12K^T\left(\nabla^2p\right)K\right)\nabla^Tp \triangle t,
\end{align*}
respectively, where $\nabla^2p$ is the Hessian matrix of $p$ and $\otimes$ is the Kronecker product. Collecting the $O(\Delta z)$ and $O(\Delta t)$ terms in $II_1-II_4$, we obtain two identities:
\begin{align}\label{eqn-K}\notag
\mathcal{O}(\Delta z):\quad0=&pK^T\nabla^2 p+p\nabla^Tp(\nabla K^T)^T+p^2\nabla^Th\\\notag
	&+p^2\left[\textup{tr}\left(\frac{\partial(\nabla K^T)}{\partial x_i}\right)\right]_{i=1,\cdots,d}\\
	&-(\nabla^Tp)K\nabla^Tp,
\end{align}
and 
\begin{align}\label{eqn-u}\notag
\mathcal{O}(\Delta t):\quad0=&pu^T\nabla^2p+\frac12pK^T\left(\nabla^T\otimes\left(\nabla^2p\right)\right)(I\otimes K)\\\notag
	&+p\nabla^Tp(\nabla u^T)^T+pK^T\nabla^2p(\nabla K^T)^T\\\notag
	&-p^2h\nabla^Th+p(\nabla^Tp)K\nabla^Th+p^2K^T\nabla^2h\\
	&+p^2\nabla^Th(\nabla K^T)^T\\\notag
	&+p(\nabla^Tp)K\left[\textup{tr}\left(\frac{\partial(\nabla K^T)}{\partial x_i}\right)\right]_{i=1,\cdots,d}\\\notag
	&+p^2\left[\textup{tr}\left(\frac{\partial(\nabla u^T)}{\partial x_i}-\nabla K^T\frac{\partial(\nabla K^T)}{\partial x_i}\right)\right]_{i=1,\cdots,d}\\\notag
	&-(\nabla^Tp)u\nabla^Tp-\frac12K^T(\nabla^2p)K\nabla^Tp.
\end{align}
\begin{proposition}\label{prop-2.2}
	The control input pair $(K,u)$ satisfies the equations:
\begin{align}\label{eqn-K-3}
	\nabla^T\left(\frac1p\nabla^T(pK)\right)=&-\nabla^Th,
\end{align}
and
\begin{align}\label{eqn-u-2}\notag
			&-\nabla^T\left(\frac1p\nabla^T(pu)\right)\\\notag
=&\nabla^T\left\{-\frac12h^2-\frac1{2p}(K^T\nabla^2pK)\right.\\\notag
	&\phantom{\nabla^Taa}+\frac12K^T\nabla(\log p)\nabla^T(\log p)K\\\notag
	&\phantom{\nabla^Taa}-K^T(\nabla K^T)\nabla(\log p)\\
	&\phantom{\nabla^Taa}\left.-\nabla^T(\nabla^TK)K-\frac12\textup{tr}\left[(\nabla K^T)(\nabla K^T)^T\right]\right\},
\end{align}
respectively.
\end{proposition}

The proof of this proposition is extremely long and involved. To avoid distraction, we postpone this tedious computations to Appendix B.

\begin{theorem}
	The {\em admissible} control input pair $(K,u)$ satisfies the equations:
\begin{align}\label{eqn-K-2}
	\nabla^T(pK)=&-(h-\hat{h})p,
\end{align}
and
\begin{align}\label{eqn-u-5}
	\nabla^T(pu)
=(h-\hat{h})\hat{h}p+\frac12\sum_{i,j=1}^d\frac{\partial^2}{\partial x_i\partial x_j}\left[p(KK^T)_{ij}\right],
\end{align}
where $\hat{h}=\int hpdx$.
\end{theorem}
\begin{IEEEproof}
Starting from \eqref{eqn-K-3}, integrating over $\R$ for each $x_i$, $i=1,\cdots,d$, yields that
\begin{align}\label{eqn-determine c}
	-\nabla^T(pK)=(h-c_i)p,
\end{align}
$i=1,\cdots,d$, where $c_i$ is a constant with respect to $x_i$. Thus, $c_i=c$, for all $i$, where $c$ is a constant to be determined. Recall that $K$ is admissible, i.e. $K\in L^2(\R^d;p)$, then $\int_{\R^d}|K|pdx\leq\left(\int_{\R^d}|K|^2pdx\right)^{1/2}<\infty$, thus, $\lim_{|x|\rightarrow\infty}pK(x)=0$. Equation \eqref{eqn-K-2} is obtained by integrating \eqref{eqn-determine c} over $\R^d$, so that
\begin{equation}\label{eqn-hath=c}
	0=\int (h-c)pdx\quad\Rightarrow\quad c=\hat{h}.
\end{equation}
To show \eqref{eqn-u-5}, integrating \eqref{eqn-u-2} once, one obtain that
\begin{align}\label{eqn-u-4}\notag
	\nabla^T(pu)=&\frac12h^2p+\frac1{2}(K^T\nabla^2pK)\\\notag
	&-\frac{p}2K^T\nabla(\log p)\nabla^T(\log p)K\\\notag
	&+pK^T(\nabla K^T)\nabla(\log p)\\\notag
	&+p\nabla^T(\nabla^TK)K\\
	&+\frac{p}2\textup{tr}\left[(\nabla K^T)(\nabla K^T)^T\right]+C_1p,
\end{align}
where $C_1$ is a constant to be determined later. Let us take a look at the third term on the right-hand side of \eqref{eqn-u-4}:
\begin{align}\label{eqn-thm1}\notag
	&-\frac p2K^T\nabla(\log p)\nabla^T(\log p)K\\\notag
\overset{\eqref{eqn-K-2}}=&-\frac p2\left[(h-\hat{h})+\nabla^TK\right]^2\\\notag
	=&-\frac p2h^2+(h-\hat{h})\hat{h}p+\frac p2\hat{h}^2+p\nabla^T(\log p)K\nabla^TK\\
	&+\frac p2(\nabla^TK)^2.
\end{align} 
Substituting \eqref{eqn-thm1} back to \eqref{eqn-u-4}, we have
\begin{align}\label{eqn-u0}\notag
	&\nabla^T(pu)\\\notag
	=&\frac1{2}(K^T\nabla^2pK)+(h-\hat{h})\hat{h}p+\frac p2\hat{h}^2+p\nabla^T(\log p)K\nabla^TK\\\notag
	&+\frac p2(\nabla^TK)^2+pK^T(\nabla K^T)\nabla(\log p)+p\nabla^T(\nabla^TK)K\\
	&+\frac{p}2\textup{tr}\left[(\nabla K^T)(\nabla K^T)^T\right]+C_1p\\\notag
=&(h-\hat{h})\hat{h}p+\frac12\sum_{i,j=1}^d\frac{\partial^2}{\partial x_i\partial x_j}\left[p(KK^T)_{ij}\right]+\frac{\hat{h}^2}2p+C_1p,
\end{align}
where the last equality follows by direct computation:
\begin{align}\label{eqn-lm2}\notag
	&\sum_{i,j=1}^d\frac{\partial^2}{\partial x_i\partial x_j}\left[p(KK^T)_{ij}\right]\\\notag
=&\sum_{i,j=1}^d\left\{\frac{\partial^2p}{\partial x_i\partial x_j}K_iK_j+2\frac{\partial p}{\partial x_i}\frac{\partial K_i}{\partial x_j}K_j+2\frac{\partial p}{\partial x_i}K_i\frac{\partial K_j}{\partial x_j}\right.\\\notag
	&\phantom{\sum_{i,j=1}^d[]}\left.+2p\frac{\partial K_j}{\partial x_i\partial x_j}K_j+p\frac{\partial K_i}{\partial x_i}\frac{\partial K_j}{\partial x_j}+p\frac{\partial K_i}{\partial x_j}\frac{\partial K_j}{\partial x_i}\right\}\\
=&K^T\nabla^2pK+2\nabla^TpK\nabla^TK+2K^T(\nabla K^T)\nabla p\\\notag
	&+p(\nabla^TK)^2+2p\nabla^T(\nabla^TK)K+p\textup{tr}\left[(\nabla K^T)(\nabla K^T)^T\right].
\end{align} 
Recall that $u$ is admissible, i.e. $u\in L^2(\R^d;p)$, then $\int_{\R^d}|u|pdx\leq\left(\int_{\R^d}|u|^2pdx\right)^{\frac12}<\infty$, thus $\lim_{|x|\rightarrow\infty}pu(x)=0$. Additionally, 
\begin{align*}
	&\int_{\R^d}\sum_{i,j}\frac{\partial^2}{\partial x_i\partial x_j}(pK_iK_j)dx\\
	=&\sum_{i,j}\int_{\R^{d-1}}\frac{\partial}{\partial x_j}\int_{\R}\frac{\partial}{\partial x_i}(pK_iK_j)dx_id(x\backslash x_j)\\
	=&\sum_{i,j}\int_{\R^{d-1}}\frac{\partial}{\partial x_j}\left[(pK_iK_j)|_{|x_i|\rightarrow\infty}\right]d(x\backslash x_j)\rightarrow0,
\end{align*}
if $||K||_\infty<\infty$. This boundedness of $K$ is followed by the admissible condition, since $K\in H^l(\R^d;p)\subset L^\infty(\R^d)$, $l\geq\lfloor\frac d2\rfloor+1$, by Sobolev embedding theorem (Theorem \ref{thm-Sobolev}). The constant $C_1=-\frac{\hat{h}^2}2$ is obtained by integrating \eqref{eqn-u0} over $\R^d$, similarly as the procedure in obtaining $c$ in \eqref{eqn-hath=c}.
\end{IEEEproof}
\begin{theorem}[Consistency]
	Suppose the admissible control input $(K,u)$ are obtained according to \eqref{eqn-K-2} and \eqref{eqn-u-5}, respectively, then provided $p(x,0)=p^*(x,0)$, we have for all $t>0$,
\begin{align*}
	p(x,t)=p^*(x,t).
\end{align*}
\end{theorem}
\begin{IEEEproof}	
Notice that the KFE of $X_t^i$ given the filtration $\mathcal{Z}_t$ is 
\begin{align*}
	dp\overset{\eqref{eqn-KFE}}=&\mathcal{L}^*pdt-\nabla^T(pK)dZ_t-\nabla^T(pu)dt\\
	&+\frac12\sum_{i,j=1}^d\frac{\partial^2}{\partial x_i\partial x_j}\left[p(KK^T)_{ij}\right]dt\\
	\overset{\eqref{eqn-K-2},\eqref{eqn-u-5}}=&\mathcal{L}^*pdt+(h-\hat{h})pdZ_t-(h-\hat{h})\hat{h}pdt\\
	=&\mathcal{L}^*pdt+(h-\hat{h})(dZ_t-\hat{h}dt)p,
\end{align*}
where $\L^*$ is defined in \eqref{eqn-L^*}, which is exactly the Kushner-Stratonovich equation \eqref{eqn-Kushner} of $X_t$. 
\end{IEEEproof}

\section{Existence and uniqueness of the control input in suitable space}

In this section, we shall discuss the existence and uniqueness of the solution in the suitable functin space to \eqref{eqn-K-2} and \eqref{eqn-u-5} under certain conditions. 

We first investigate the weak solution to \eqref{eqn-K-2} with $K=\nabla\phi$ such that
\begin{equation}
	\int\nabla^T\psi\nabla\phi pdx=\int(h-\hat{h})\psi pdx,
\end{equation}
for all $\psi\in H_0^1(\R^d;p)$, which has the norm defined 
\begin{align*}
	||\psi||_{H_0^1(\R^d;p)}^2=||\psi||_{L^2(\R^d;p)}^2+||\nabla\psi||_{L^2(\R^d;p)}^2,
\end{align*}
with $\lim_{|x|\rightarrow\infty}\psi(x)=0$.

Assume that
\begin{enumerate}
	\item[(As-1)] $h\in H^{k}(\R^d;p)(:=W^{k,2}(\R^d;p))$, for some $k\geq0$; 
	\item[(As-2)] Poincar\'{e}-type inequality: there exists a constant $C>0$ such that for any $\phi\in H_0^1(\R^d;p)$,
\begin{equation}\label{eqn-Poincare2}
	||\phi||_{L^2(\R^d;p)}\leq C||\nabla\phi||_{L^2(\R^d;p)}
\end{equation}
	\item[(As-3)] $\nabla^l(\log p)\in L^\infty(\R^d)$, for $2\leq l\leq k+1$. Here, $k$ is the one in (As-1), with the convention that $l\in\emptyset$, if $k<1$.
\end{enumerate}

In \cite{YLMM}, Yang et. al. gave the similar assumptions as above to guarantee the existence and uniqueness of the solution $\phi$ to \eqref{eqn-K-2}, but they have no discussion on the uniqueness of the control input $u$. They gave the explicit expression \eqref{eqn-explicit solution u} of $u$ and verified that $u\in L^1(\R^d;p)$. 

We point out that Assumption (As-2) may not hold in $\R^d$. This has been proven in Lemma 10.2(ii), \cite{T} with the unweighted norm. In the following lemma, we show the similar result in the weighted Sobolev space $W_0^{1,q}(\Omega;p)$, $1\leq q<\infty$:
\begin{lemma}
	If for some small $\epsilon\ll1$, $||\nabla(\log p)||_\infty\leq q(1-\epsilon)$ and $p(x)>0$ in $\Omega$, then Poincar\'{e}-type inequality
\begin{equation}\label{eqn-Poincare}
	||u||_{L^q(\Omega;p)}\leq C||\nabla u||_{L^q(\Omega;p)}
\end{equation}
does not hold for $u\in W^{1,q}_0(\Omega;p)$, if $\Omega$ contains arbitrarily large balls, i.e., if there exists a sequence $r_n\rightarrow\infty$ and points $x_n\rightarrow\Omega$ such that the ball centered at $x_n$ with radius $r_n$ is in $\Omega$, i.e., $B(x_n,r_n)\in\Omega$. Here, $W_0^{1,q}(\Omega;p)$ is the Sobolev space with 
\begin{equation}
	||\circ||^q_{W^{1,q}_0(\Omega;p)}=||\circ||^q_{L^q(\Omega;p)}+||\nabla\circ||^q_{L^q(\Omega;p)}
\end{equation}
and $\circ(x)\equiv0$ on $\partial\Omega$.
\end{lemma}
\begin{IEEEproof}(By contradiction)
	Let $\gamma(x)\in C_c^\infty(\R^d)$ with $\gamma\not\equiv0$ and $\textup{support}(\gamma)\subset B(0,1)$, then one defines $u_n(x)=\gamma(\frac{x-x_n}{r_n})\frac1{p^{\frac1q}(x)}$ is also compactly supported in $\Omega$. Thus, one has 
\begin{align}\label{eqn-unLq}\notag
	||u_n||_{L^q(\Omega;p)}=&\left[\int_\Omega\left|\gamma\left(\frac{x-x_n}{r_n}\right)\right|^qdx\right]^{1/q}\\\notag
	=&r_n^{d/q}\left[\int_{B(0,1)}|\gamma(y)|^qdy\right]^{1/q}\\
	=&r_n^{d/q}||\gamma||_{L^q(B(0,1))},
\end{align}
while
\begin{align}\label{eqn-nablaunLq}\notag
	&||\nabla u_n||_{L^q(\Omega;p)}\\\notag
\leq&\frac1{r_n}\left[\int_\Omega\left|\nabla\gamma\left(\frac{x-x_n}{r_n}\right)\right|^qdx\right]^{1/q}\\\notag
	&+\frac1q\left[\int_\Omega\left|\gamma\left(\frac{x-x_n}{r_n}\right)\right|^q|\nabla(\log p)|^qdx\right]^{1/q}\\\notag
	=&r_n^{-1+\frac dq}\left[\int_{B(0,1)}|\nabla\gamma(y)|^qdy\right]^{\frac1q}\\\notag
	&+\frac{r_n^{\frac dq}}q\left[\int_{B(0,1)}|\gamma(y)|^q|\nabla(\log p)|^qdy\right]^{\frac1q}\\\notag
	\leq&r_n^{-1+\frac dq}||\nabla\gamma||_{L^q(B(0,1))}+\frac{r_n^{\frac dq}}q||\nabla(\log p)||_\infty||\gamma||_{L^q(B(0,1))}\\
	\leq&r_n^{-1+\frac dq}||\nabla\gamma||_{L^q(B(0,1))}+r_n^{\frac dq}(1-\epsilon)||\gamma||_{L^q(B(0,1))},
\end{align}
where the last inequality follows from the assumption on $||\nabla(\log p)||_\infty$. Suppose Poincar\'{e}-type inequality \eqref{eqn-Poincare} holds for all $u_n\in W^{1,q}_0(\Omega;p)$, then there exists a constant $C$ independent of $n$ such that 
\begin{align*}
	0<||\gamma||_{L^q(B(0,1))}\leq\frac C{r_n\epsilon}||\nabla\gamma||_{L^q(B(0,1))}\rightarrow0,
\end{align*}
as  $r_n\rightarrow\infty$. The contradiction is arrived.
\end{IEEEproof}

However, it has been verified that Assumption (As-2) does hold for the nonlinear non-Gaussian case with a constant signal model in \cite{LMMR}. It is an interesting question that under what conditions, the Poincar\'{e}-type inequality is guaranteened for $W^{1,q}_0(\R^d;p)$. This is our future project.

\subsection{The existence and uniqueness of $\phi\in H_0^{k+1}(\R^d;p)$}

\begin{theorem}\label{thm-about k}
	Under Assumptions (As-1)-(As-3), equation \eqref{eqn-K-2} has a unique weak solution $\phi\in H_0^{k+1}(\R^d;p)$, where $k$ is the one in (As-1).
\end{theorem}

It follows easily from Lax-Milgram theorem that there exists a unique solution $\phi\in H_0^1(\R^d;p)$. Or one can consult the proof of Theorem 2, section A.2, \cite{YLMM}. We shall show that $\phi\in H_0^{k+1}(\R^d;p)$, $k\geq1$, in Lemma \ref{lem-apriori esstimate} later. But before that, we need the following lemma:
\begin{lemma}\label{lem-1}
	If $\phi\in C^\infty(\R^d)$ satisfies \eqref{eqn-K-2} with $K=\nabla\phi$, we have 
\begin{align}\label{eqn-DkeqnK}
	-\nabla^T\left[p\nabla\left(\nabla_{\sigma_l}^l\phi\right)\right]
	=G_{\sigma_l}^lp,
\end{align}
for $l\geq1$, where $\nabla_{\sigma_l}^l\circ$ is defined in \eqref{eqn-nablasigmak}, $\sigma_l\in \{1,\cdots,d\}^l$ and 
\begin{align}\label{eqn-G1}
	G_{\sigma_1}^1:=&\nabla^T\left[\frac{\partial(\log p)}{\partial x_{\sigma_1(1)}}\right]\nabla\phi+\frac{\partial h}{\partial x_{\sigma_1(1)}},\\\label{eqn-Gk}
	G_{\sigma_l}^l:=&\nabla^T\left[\frac{\partial(\log p)}{\partial x_{\sigma_l(l)}}\right]\nabla\left(\nabla_{\sigma_l|_{1:l-1}}^{l-1}\phi\right)+\frac{\partial G^{l-1}_{\sigma_l|_{1:l-1}}}{\partial x_{\sigma_l(l)}},	
\end{align}
for $l\geq2$, with $\sigma_l|_{1:l-1}=(\sigma_l(1),\cdots,\sigma_l(l-1))\in\{1,\cdots,d\}^{l-1}$ is the first $l-1$ component of $\sigma_l$, where
\begin{align}\label{eqn-parialGk}\notag
	&\frac{\partial G^{l-1}_{\sigma_l|_{1:l-1}}}{\partial x_{\sigma_l(l)}}\\\notag
=&\sum_{i=2}^{l}\nabla^{l-i+1}_{\sigma_l|_{i:l}}\left[\nabla^T\left(\frac{\partial(\log p)}{\partial x_{\sigma_{l}(i-1)}}\right)\nabla\left(\nabla^{i-2}_{\sigma_{l}|_{1:i-2}}\phi\right)\right]\\
	&+\nabla^l_{\sigma_{l}}h,
\end{align}
with the convention that $\nabla^0_\cdot=Id$, where $Id$ represents the identity mapping.
\end{lemma}

To avoid the distraction, we append the proof of this lemma in Appendix D.

\begin{lemma}\label{lem-apriori esstimate}
	Under Assumption (As-1)-(As-3), if the weak solution of \eqref{eqn-K-2} with $K=\nabla\phi$, $\phi\in H_0^{1}(\R^d;p)$, then $\phi$ is actually is in $H^{k+1}_0(\R^d;p)$, for $k$ in (As-1). Moreover, it  satisfies that for any $l\geq1$ and $\sigma_l\in\{1,\cdots,d\}^l$ such that
\begin{align}\label{eqn-estimate of dk+1L2norm}\notag
	&||\nabla^{l+1}\phi||^2_{L^2(\R^d;p)}\\
	 \leq&||\nabla^l\phi||^2_{L^2(\R^d;p)}+\sum_{\sigma_l\in\{1,\cdots,d\}^l}||G_{\sigma_l}^l||^2_{L^2(\R^d;p)},
\end{align}
with
\begin{align}\label{eqn-||G^k||}\notag
	&\left|\left|G_{\sigma_l}^l\right|\right|^2_{L^2(\R^d;p)}\\\notag
\lesssim&\sum_{i=2}^l\sum_{m=0}^{l-i+1}\sum_{\sigma_m\subset\sigma_l}\left|\left|\nabla^T\left[\nabla_{(\sigma_m,\sigma_l(i-1))}^{m+1}(\log p)\right]\right|\right|_\infty\\\notag
	&\phantom{\sum_{i=1}^l\sum_{m=0}^{l-i+1}\sum_{\sigma_m\subset\sigma_l}aa}\cdot\left|\left|\nabla\left[\nabla^{l-m-1}_{(\sigma_l|_{1:i-2},\sigma_l|_{i:l}\backslash\sigma_m)}\phi\right]\right|\right|^2_{L^2(\R^d;p)}\\
	&+||\nabla^l_{\sigma_{l}}h||^2_{L^2(\R^d;p)},
\end{align}
where
\begin{align}\label{eqn-def of L2}
	||\nabla^i\phi||^2_{L^2(\R^d;p)}=\sum_{\sigma_i\in\{1,\cdots,d\}^i}\int|\nabla^i_{\sigma_i}\phi|^2pdx.
\end{align}
\end{lemma}
\begin{IEEEproof}
	We first claim that for $l\geq1$, $\sigma_l\in\{1,\cdots,d\}^l$,
\begin{align}\label{eqn-control of dk+1}
	\int|\nabla(\nabla^l_{\sigma_l}\phi)|^2pdx\leq\int|\nabla_{\sigma_l}^l\phi||G_{\sigma_l}^l|pdx,
\end{align}
holds if $\phi\in H_0^{l+1}(\R^d;p)$. 

Let $\beta(x)\geq0$ be a smooth, compactly supported, radially decreasing function in $\R^d$, with $\beta(0)=1$. Let $s>0$. For any $l\geq1$ and $\sigma_l\in\{1,\cdots,d\}^l$, multiplying \eqref{eqn-DkeqnK} with $\beta^2(sx)\nabla^l_{\sigma_l}\phi(x)$ and integrating it over $\R^d$, yields that
\begin{align}\label{eqn-intbeta}\notag
	&\int\beta^2(sx)\nabla_{\sigma_l}^l\phi G_{\sigma_l}^lpdx\\\notag
=&\int\nabla^T\left[\beta^2(sx)\nabla_{\sigma_l}^l\phi\right]\nabla(\nabla_{\sigma_l}^l\phi)pdx\\\notag
=&\int2s\beta(sx)\nabla^T\beta(sx)\nabla_{\sigma_l}^l\phi\nabla(\nabla_{\sigma_l}^l\phi)pdx\\
	&+\int\beta^2(sx)\left|\nabla(\nabla_{\sigma_l}^l\phi)\right|^2pdx,
\end{align}
where the first equality follows by integration by parts, and $\nabla^l_{\sigma_l}$ is defined as \eqref{eqn-nablasigmak}. The first term on the right-hand side of \eqref{eqn-intbeta} can be controlled as
\begin{align}\label{eqn-control1}\notag
	&\left|\int2s\beta(sx)\nabla^T\beta(sx)(\nabla_{\sigma_l}^l\phi)\nabla(\nabla_{\sigma_l}^l\phi)pdx\right|\\
	\leq&2s||\nabla^T\beta||_\infty\int\beta(sx)\left|\nabla_{\sigma_l}^l\phi\right|\left|\nabla(\nabla_{\sigma_l}^l\phi)\right|pdx\\\notag
	\leq&s||\nabla^T\beta||_\infty\left(\int|\nabla_{\sigma_l}^l\phi|^2pdx+\int\beta^2(sx)|\nabla(\nabla_{\sigma_l}^l\phi)|^2pdx\right),
\end{align}
where the last inequality follows by Cauchy-Schwarz inequality. Substituting \eqref{eqn-control1} back to \eqref{eqn-intbeta}, we obtain that
\begin{align*}
	&\int|G_{\sigma_l}^l||\nabla_{\sigma_l}^l\phi|pdx\\
	\geq&\int\beta^2(sx)G_{\sigma_l}^l(\nabla_{\sigma_l}^l\phi)pdx\\
\overset{\eqref{eqn-intbeta}}\geq&(1-s||\nabla^T\beta||_\infty)\int\beta^2(sx)\left|\nabla(\nabla_{\sigma_l}^l\phi)\right|^2pdx\\
	&-s||\nabla^T\beta||_\infty\int|\nabla_{\sigma_l}^l\phi|^2pdx\\
\geq&\int\left|\nabla(\nabla_{\sigma_l}^l\phi)\right|^2pdx\\
	&-s||\nabla^T\beta||_\infty\left[\int\left|\nabla(\nabla_{\sigma_l}^l\phi)\right|^2pdx+\int|\nabla_{\sigma_l}^l\phi|^2pdx\right],
\end{align*}
where the last inequality is due to $\beta(0)=1$ and $\beta$ is radially decreasing. Thus, \eqref{eqn-control of dk+1} follows immediately by letting $s\rightarrow0$ and dominated convergence theorem.

In the sequel, we show that \eqref{eqn-estimate of dk+1L2norm} holds, for $1\leq l\leq k$, here $k$ is the one in (As-1). In fact, the direct computation yields that
\begin{align*}
	&||\nabla^{l+1}\phi||^2_{L^2(\R^d;p)}\\
	=&\sum_{{\sigma_{l+1}}\in\{1,\cdots,d\}^{l+1}}\int|\nabla^{l+1}_{\sigma_{l+1}}\phi|^2pdx\\
	=&\sum_{{\sigma_{l+1}}|_{1:l}\in\{1,\cdots,d\}^l}\int|\nabla(\nabla^l_{{\sigma_{l+1}}|_{1:l}}\phi)|^2pdx\\
	\overset{\eqref{eqn-control of dk+1}}\leq&\sum_{{\sigma_{l+1}}|_{1:l}\in\{1,\cdots,d\}^l}\left[\int|\nabla_{{\sigma_{l+1}}|_{1:l}}^l\phi|^2pdx\right]^{\frac12}\\
	&\phantom{\sum_{{\sigma_{l+1}}|_{1:l}\in\{1,\cdots,d\}^l}aa}\cdot\left[\int|G_{{\sigma_{l+1}}|_{1:l}}^l|^2pdx\right]^{\frac12}\\
	\lesssim&\sum_{{\sigma_{l+1}}|_{1:l}\in\{1,\cdots,d\}^l}\int|\nabla_{{\sigma_{l+1}}|_{1:l}}^l\phi|^2pdx\\
	&+\sum_{{\sigma_{l+1}}|_{1:l}\in\{1,\cdots,d\}^l}\int|G_{{\sigma_{l+1}}|_{1:l}}^l|^2pdx,
\end{align*}
where the last two inequalities follow from H\"{o}lder's inequality and Cauchy-Schwarz inequality, respectively. Thus, the right-hand side of \eqref{eqn-estimate of dk+1L2norm} is obtained by recalling the definition of $L^2$-norm \eqref{eqn-def of L2}.

From \eqref{eqn-estimate of dk+1L2norm}, we claim by induction that $\nabla^{l+1}\phi, G_{\sigma_{l+1}}^{l+1}\in L^2(\R^d;p)$, for any $l\geq1$, provided $G_{{\sigma_1}(1)}^1,\nabla\phi\in L^2(\R^d;p)$. 

We begin with examining $||\nabla\phi||_{L^2(\R^d;p)}<\infty$ and $||G_{\sigma_1(1)}^1||_{L^2(\R^d;p)}<\infty$. Multiplying \eqref{eqn-K-2} with $K=\nabla\phi$ by $\phi$ and integrating over $\R^d$:
\begin{align*}
	||\nabla\phi||_{L^2(\R^d;p)}^2=&-\int\nabla^T(p\nabla\phi)\phi dx
\overset{\eqref{eqn-K-2}}=\int(h-\hat{h})\phi pdx\\
	\leq&\left(||h||_{L^2(\R^d;p)}+|\hat{h}|\right)||\phi||_{L^2(\R^d;p)}\\
	\overset{\eqref{eqn-Poincare2}}\leq&C\left(||h||_{L^2(\R^d;p)}+|\hat{h}|\right)||\nabla\phi||_{L^2(\R^d;p)}.
\end{align*}
Therefore, we obtain that $||\nabla\phi||_{L^2(\R^d;p)}<\infty$, provided $h\in L^2(\R^d;p)$.
%\begin{align}\label{eqn-estimate nablaphi}
	%||\nabla\phi||_{L^2(\R^d;p)}\leq\frac1\lambda\left(||h||_{L^2(\R^d;p)}+\hat{h}\right)<\infty.
%\end{align}
Now, we check whether $G_{\sigma_1(1)}^1\in L^2(\R^d;p)$:
\begin{align}\label{eqn-estimate G11}
	&||G_{\sigma_1(1)}^1||_{L^2(\R^d;p)}^2\\\notag
	\overset{\eqref{eqn-G1}}\leq&\left|\left|\nabla^T\left[\frac{\partial(\log p)}{\partial x_{\sigma_1(1)}}\right]\nabla\phi\right|\right|^2_{L^2(\R^d;p)}+\left|\left|\frac{\partial h}{\partial x_{\sigma_1(1)}}\right|\right|_{L^2(\R^d;p)}^2\\\notag
	\overset{\textup{(As-3)}}\leq&\left|\left|\nabla^2(\log p)\right|\right|_\infty||\nabla\phi||^2_{L^2(\R^d;p)}+\left|\left|\nabla h\right|\right|_{L^2(\R^d;p)}^2<\infty,
\end{align} 
provided that $\nabla\phi,\nabla h\in L^2(\R^d;p)$ and (As-3).

By induction, suppose that for all $\sigma_{l-1}\in \{1,\cdots,d\}^{l-1}$, $G_{\sigma_{l-1}}^{l-1}, \nabla^{l-1}\phi\in L^2(\R^d;p)$, then we have $\nabla^l\phi\in L^2(\R^d;p)$, $l\geq2$, by \eqref{eqn-estimate of dk+1L2norm}, and
\begin{align}\label{eqn-||Gl||2}\notag
		&\left|\left|G_{\sigma_l}^l\right|\right|^2_{L^2(\R^d;p)}\\\notag
\overset{\eqref{eqn-Gk}}\leq&\left|\left|\nabla^T\left[\frac{\partial(\log p)}{\partial x_{\sigma_l(l)}}\right]\right|\right|_\infty\left|\left|\nabla\left(\nabla_{\sigma_l|_{1:l-1}}^{l-1}\phi\right)\right|\right|^2_{L^2(\R^d;p)}\\
	&+\left|\left|\frac{\partial G^{l-1}_{\sigma_l|_{1:l-1}}}{\partial x_{\sigma_l(l)}}\right|\right|_{L^2(\R^d;p)}^2,
\end{align}
since
\begin{align}\label{eqn-partial Gl}\notag
	&\left|\left|\frac{\partial G_{\sigma_l|_{1:l-1}}^{l-1}}{\partial x_{\sigma_l(l)}}\right|\right|^2_{L^2(\R^d;p)}\\\notag
\overset{\eqref{eqn-parialGk}}\lesssim&\sum_{i=2}^{l}\left|\left|\nabla^{l-i+1}_{\sigma_l|_{i:l}}\left[\nabla^T\left(\frac{\partial(\log p)}{\partial x_{\sigma_{l}{(i-1)}}}\right)\nabla\left(\nabla^{i-2}_{\sigma_{l}|_{1:i-2}}\phi\right)\right]\right|\right|^2_{L^2(\R^d;p)}\\\notag
	&+||\nabla^l_{\sigma_{l}}h||^2_{L^2(\R^d;p)}\\\notag
	\lesssim&\sum_{i=2}^l\sum_{m=0}^{l-i+1}\sum_{\sigma_m\subset\sigma_l}\left|\left|\nabla^T\left[\nabla_{(\sigma_m,\sigma_l(i-1))}^{m+1}(\log p)\right]\right|\right|_\infty\\\notag
	&\phantom{\sum_{i=2}^l\sum_{m=0}^{l-i+1}\sum_{\sigma_m\subset\sigma_l}aa}\cdot\left|\left|\nabla\left[\nabla^{l-m-1}_{(\sigma_l|_{1:i-2},\sigma_l|_{i:l}\backslash\sigma_m)}\phi\right]\right|\right|^2_{L^2(\R^d;p)}\\
	&+||\nabla^l_{\sigma_{l}}h||^2_{L^2(\R^d;p)},
\end{align}
where $\sigma_{m}$ is the sub-vector of $\sigma$ of size $m$, while $\sigma\backslash\sigma_{m}$ is the rest of the vector $\sigma$ after removing $\sigma_{m}$, with the convention that $\sigma_0=\emptyset$. Therefore, \eqref{eqn-||G^k||} follows by substituting \eqref{eqn-partial Gl} back to \eqref{eqn-||Gl||2}. Consequently, $G_{\sigma_l}^l\in L^2(\R^d;p)$, provided $\nabla^m\phi,\,\nabla ^mh\in L^2(\R^d;p)$ and $||\nabla^m(\log p)||_\infty<\infty$, for all $m\leq l$.

Under Assumption (As-1) and (As-3), we conclude that $\phi\in H_0^{k+1}(\R^d;p)$, if $\phi\in H_0^1(\R^d;p)$. Here, $k$ is the one in (As-1).
\end{IEEEproof}

\subsection{The existence and uniqueness of the solution $\varphi\in H^1_0(\R^d;p)$}

In \cite{YLMM}, the expression of the solution $u$ is directly given without any derivation. There is NO uniqueness result of this solution $u$. In this subsection, we confirm that the explicit solution given in \cite{YLMM} is the unique one in the space $L^2(\R^d;p)$ under Assumption (As-1)-(As-3), with $k\geq\lfloor\frac{d+2}4\rfloor$ in (As-1).

Similarly as we did for \eqref{eqn-K-2}, we seek the weak solution to \eqref{eqn-u-5} with $u=\nabla\varphi$ such that
\begin{align}\label{eqn-weak varphi}
\int\nabla^T\varphi\nabla\psi pdx
	=&\hat{h}\int h\psi pdx-\hat{h}^2\int\psi pdx\\\notag
	&+\int\frac12\sum_{i,j=1}^d\frac{\partial^2}{\partial x_i\partial x_j}\left(pK_iK_j\right)\psi dx,
\end{align}
for $\psi\in H_0^1(\R^d;p)$.

\begin{theorem}\label{thm-3.2}
	Under Assumption (As-1)-(As-3) with $k\geq\lfloor\frac{d+2}4\rfloor$ in (As-1), equation \eqref{eqn-u-5} has a unique solution $u=\nabla\varphi$, with $\varphi\in H^1_0(\R^d)$ .
\end{theorem}
\begin{IEEEproof}
	The existence and uniqueness of $\varphi\in H^{1}_0(\R^d;p)$ is guaranteened by Riesz represetation theorem. In fact, the inner product of $(\varphi,\psi)$ is defined as the integral on the left-hand side of \eqref{eqn-weak varphi}. We only need to check its boundedness. Let us first look at the last term on the right-hand side of \eqref{eqn-weak varphi}:
\begin{align}\label{eqn-V3}\notag
	&\int\sum_{i,j=1}^d\frac{\partial^2}{\partial x_i\partial x_j}\left(pK_iK_j\right)\psi dx\\\notag
=&-\int\sum_{i,j=1}^d\frac{\partial}{\partial x_i}(pK_iK_j)\frac{\partial  \psi}{\partial x_j}dx\\\notag
	=&-\int\sum_{i,j=1}^d\frac{\partial}{\partial x_i}(pK_i)K_j\frac{\partial\psi}{\partial x_j}dx
-\int\sum_{i,j=1}^dpK_i\frac{\partial K_j}{\partial x_i}\frac{\partial \psi}{\partial x_j}dx\\\notag
	=&-\int\nabla^T(p\nabla\phi)\nabla^T\phi\nabla\psi dx
-\int \nabla^T\phi\nabla^2\phi\nabla\psi pdx\\
	\overset{\eqref{eqn-K-2}}=&\int(h-\hat{h})\nabla^T\phi\nabla\psi pdx-\int \nabla^T\phi\nabla^2\phi\nabla\psi pdx.
\end{align}
By Theorem \ref{thm-about k}, we have $\phi\in H^{k+1}_0(\R^d;p)$. Thus, the inner product can be bounded
\begin{align}\label{eqn-estimate inner product}
	&|\langle\varphi,\psi\rangle|\\\notag
\overset{\eqref{eqn-weak varphi},\eqref{eqn-V3}}\lesssim&|\hat{h}|\,||h||_{L^2(\R^d;p)}||\psi||_{L^2(\R^d;p)}+|\hat{h}|^2||\psi||_{L^2(\R^d;p)}\\\notag
	&+||h||_{L^r(\R^d;p)}||\nabla\phi||_{L^{r^*}(\R^d;p)}||\nabla\psi||_{L^2(\R^d;p)}\\\notag
	&+||\nabla\phi||_{L^s(\R^d;p)}||\nabla^2\phi||_{L^{s^*}(\R^d;p)}||\nabla\psi||_{L^2(\R^d;p)},
\end{align}
where $\frac1r+\frac1{r^*}=\frac1s+\frac1{s^*}=\frac12$. By Sobolev embedding theorem (Theorem \ref{thm-Sobolev}), it is easy to deduce that for the second term on the right-hand side to be bounded:
\begin{align*}&
	\left.\begin{aligned}
	||h||_{L^r(\R^d;p)}:&\ \frac1r\geq\frac12-\frac kd\\	
	||\nabla\phi||_{L^{r^*}(\R^d;p)}:&\ \frac1{r^*}-\frac1d\geq\frac12-\frac {k+1}d
\end{aligned}\right\}\\
\Rightarrow&\frac12=\frac1r+\frac1{r^*}\geq1-\frac{2k}d,
\end{align*}
thus, $k\geq\frac{d}4$, and for the third term:
\begin{align*}
	&\left.\begin{aligned}
	||\nabla\phi||_{L^s(\R^d;p)}:&\ \frac1s-\frac1d\geq\frac12-\frac {k+1}d\\
	||\nabla^2\phi||_{L^{s^*}(\R^d;p)}:&\ \frac1{s^*}-\frac2d\geq\frac12-\frac {k+1}kd
\end{aligned}\right\}\\
&\Rightarrow\frac12=\frac1s+\frac1{s^*}\geq1-\frac{2k-1}d,
\end{align*}
thus, $k\geq\frac{d+2}4$. Therefore, for $k\geq\lfloor\frac{d+2}4\rfloor$, the right-hand side of \eqref{eqn-estimate inner product} is bounded. That is, 
\begin{align*}
	&|\langle\varphi,\psi\rangle|\\
\lesssim&\left[C\left(|\hat{h}|\,||h||_{L^2(\R^d;p)}+|\hat{h}|^2\right)\right.
	+||h||_{L^r(\R^d;p)}||\nabla\phi||_{L^{r^*}(\R^d;p)}\\
	&\phantom{aa}\left.+||\nabla\phi||_{L^s(\R^d;p)}||\nabla^2\phi||_{L^{s^*}(\R^d;p)}\right]||\nabla\psi||_{L^2(\R^d;p)}\\
	\leq&C||\nabla\psi||_{L^2(\R^d;p)},
\end{align*}
where the first inequality follows from Assumption (As-2). By Riesz representation theorem, there exists a unique solution $\varphi\in H_0^{1}(\R^d;p)$ such that \eqref{eqn-weak varphi} holds for all $\psi\in H^{1}_0(\R^d;p)$.
\end{IEEEproof}
\begin{remark}\label{rm-1}
	Suppose $\phi\in H_0^{k+1}(\R^d;p)$, $k\geq\lfloor\frac{d+2}4\rfloor$, and $\varphi\in H_0^1(\R^d;p)$ are the solutions in Theorem \ref{thm-about k} and  \ref{thm-3.2}, respectively, then $K\in H_0^{k}(\R^d;p)$ and $u\in L^2(\R^d;p)$ are also unique. 
\end{remark}
\begin{corollary}\label{coro-1}
	Under Assumption (As-1)-(As-3), with $k\geq\lfloor\frac{d+2}4\rfloor$ the explicit solution 
\begin{equation}\label{eqn-explicit solution u}
	u=-\frac12K(h+\hat{h})+\Omega,
\end{equation}
where $\Omega=\frac12(\nabla K^T)^TK$, given in \cite{YLMM} is indeed the unique solution in $L^2(\R^d;p)$.  
\end{corollary}
\begin{IEEEproof}
	The conclusion is arriven by examining that $u$ in \eqref{eqn-explicit solution u} is in $L^2(\R^d;p)$:
\begin{align*}
	&||u||^2_{L^2(\R^d;p)}\\
\overset{\eqref{eqn-explicit solution u}}\lesssim&\int|\nabla\phi|^2|h+\hat{h}|^2pdx+\int|\nabla^2\phi|^2|\nabla\phi|^2pdx\\
	\leq&|||h||^2_{L^{2r}(\R^d;p)}|\nabla\phi||^2_{L^{2r^*}(\R^d;p)}+|\hat{h}|^2||\nabla\phi||^2_{L^2(\R^d;p)}\\
	&+||\nabla\phi||^2_{L^{2s}(\R^d;p)}||\nabla^2\phi||_{L^{2s^*}(\R^d;p)}^2<\infty,
\end{align*}
where $\frac1r+\frac1{r^*}=\frac1s+\frac1{s^*}=1$, as argued in \eqref{eqn-estimate inner product} in the proof of Theorem \ref{thm-3.2}. As mentioned in Remark \ref{rm-1}, the solution $u\in L^2(\R^d;p)$ is unique, then \eqref{eqn-explicit solution u} is the one.
\end{IEEEproof}

In fact, if $h$ has even higher regularity, then the control input $u$ can be smoother. 
\begin{theorem}
Under Assumption (As-1) and (As-3) with $k\geq\lfloor\frac{d+2}4\rfloor$, the unique solution $u$ with expression \eqref{eqn-explicit solution u} is actually in $H^{l}_0(\R^d;p)$, with $l=\min\left\{k,\lfloor2k-1-\frac d2\rfloor\right\}$. 
\end{theorem}
\begin{IEEEproof}
	With the similar argument in the proof of Theorem \ref{thm-3.2}, we have $\phi\in H^{k+1}_0(\R^d;p)$, $k\geq\lfloor\frac{d+2}4\rfloor$. By Corollary \ref{coro-1}, $u$ in \eqref{eqn-explicit solution u} is the unique solution in $L^2(\R^d;p)$. In the sequel, we shall show that actually this unique solution $u\in H^{l}_0(\R^d;p)$, for some $l\geq1$. For any $\sigma\in\{1,\cdots,d\}^l$,
\begin{align}\label{eqn-nablal}\notag
	\nabla_\sigma^lu
	=&-\frac12\sum_{l_1=0}^l\sum_{\sigma_{l_1}\subset\sigma}\nabla(\nabla^{l_1}_{\sigma_{l_1}}\phi)\nabla^{l-l_1}_{\sigma\backslash\sigma_{l_1}}h
	-\frac12\hat{h}\nabla(\nabla^l_\sigma\phi)\\
	&+\frac12\sum_{l_1=0}^l\sum_{\sigma_{l_1}\subset\sigma}\nabla^2(\nabla^{l_1}_{\sigma_{l_1}}\phi)\nabla(\nabla^{l-l_1}_{\sigma\backslash\sigma_{l_1}}\phi),
\end{align}
where $\sigma_{l_1}$ is the sub-vector of $\sigma$ of size $l_1$ and $\sigma\backslash\sigma_{l_1}$ is the same notation as before. Taking $L^2(\R^d;p)$ norm on both sides of \eqref{eqn-nablal}, we have
\begin{align*}
	&||\nabla^l_{\sigma}u||^2_{L^2(\R^d;p)}\\
\lesssim&\sum_{l_1=0}^l\sum_{\sigma_{l_1}\subset\sigma}||\nabla(\nabla^{l_1}_{\sigma_{l_1}}\phi)||_{L^{2r_{\sigma_{l_1}}}(\R^d;p)}^2||\nabla^{l-l_1}_{\sigma\backslash\sigma_{l_1}}h||^2_{L^{2r^*_{\sigma_{l_1}}}(\R^d;p)}\\
	&+|\hat{h}|\,||\nabla(\nabla_\sigma^l\phi)||_{L^2(\R^d;p)}^2\\
	&+\sum_{l_1=0}^l\sum_{\sigma_{l_1}\subset\sigma}||\nabla^2(\nabla^{l_1}_{\sigma_{l_1}}\phi)||^2_{L^{2s_{\sigma_{l_1}}}(\R^d;p)}\\
	&\phantom{+\sum_{l_1=0}^l\sum_{\sigma_{l_1}\subset\sigma}}\cdot||\nabla(\nabla^{l-l_1}_{\sigma\backslash\sigma_{l_1}}\phi)||^2_{L^{2s^*_{\sigma_{l_1}}}(\R^d;p)}\\
	=:&V_1+V_2+V_3,
\end{align*}
where $\frac1{r_{\sigma_{l_1}}}+\frac1{r^*_{\sigma_{l_1}}}=\frac1{s_{\sigma_{l_1}}}+\frac1{s^*_{\sigma_{l_1}}}=1$, for all $l_1=0,\cdots,l$, $\sigma_{l_1}\subset\sigma\in\{1,\cdots,d\}^l$.

Next, we check that how large $l$ could be so that $V_1-V_3$ are bounded. For $V_1$:
\begin{align*}
	&\left.\begin{aligned}
	||\nabla(\nabla^{l_1}_{\sigma_{l_1}}\phi)||_{L^{2r_{\sigma_{l_1}}}(\R^d;p)}:&\ \frac1{2r_{\sigma_{l_1}}}-\frac{l_1+1}d\geq\frac12-\frac {k+1}d\\
	||\nabla^{l-l_1}_{\sigma\backslash\sigma_{l_1}}h||_{L^{2r^*_{\sigma_{l_1}}}(\R^d;p)}:&\ \frac1{2r^*_{\sigma_{l_1}}}-\frac{l-l_1}d\geq\frac12-\frac kd
\end{aligned}\right\}\\
	&\Rightarrow1=\frac1{r_{\sigma_{l_1}}}+\frac1{r^*_{\sigma_{l_1}}}\geq2-\frac{4k-2l}d\\
	&\Rightarrow l\leq2k-\frac d2.
\end{align*}
For $V_2$, it is clear to see that $l\leq k$. For $V_3$:
\begin{align*}
	&\left.\begin{aligned}
	||\nabla^2(\nabla^{l_1}_{\sigma_{l_1}}\phi)||_{L^{2s_{\sigma_{l_1}}}(\R^d;p)}:&\ \frac1{2s_{\sigma_{l_1}}}-\frac{l_1+2}d\geq\frac12-\frac {k+1}d\\
	||\nabla(\nabla^{l-l_1}_{\sigma\backslash\sigma_{l_1}}\phi)||_{L^{2s^*_{\sigma_{l_1}}}(\R^d;p)}:&\ \frac1{2s^*_{\sigma_{l_1}}}-\frac{l-l_1+1}d\geq\frac12-\frac {k+1}d
\end{aligned}\right\}\\
	&\Rightarrow1=\frac1{s_{\sigma_{l_1}}}+\frac1{s^*_{\sigma_{l_1}}}\geq2-\frac{4k-2l-2}d\\
	&\Rightarrow l\leq2k-\frac d2-1.
\end{align*}
Therefore, $l=\min\left\{k,\lfloor2k-1-\frac d2\rfloor\right\}$.
\end{IEEEproof}

\section{Conclusion}

In this paper, we re-investigated the multivariate feedback particle filter. We first show that the control input can be obtained by any f-divergence, not just the Kullback-Leibler divergence. Next, we carefully derived the equations that the control input satisfies, which has not been done in \cite{YLMM}. The derivation is for $d\geq1$ and $m=1$, but can be extended to the most general case $m\geq1$, with more involved notations and computations. We re-defined admissible for the control input, so that it can be shown that with this new definition the admissible control input exists and is unique in appropriate function spaces. Furthermore, we show that the explicit expression given in \cite{YLMM} for the control input is actually the only admissible one.

\section*{Appendix}

\setcounter{theorem}{0}
\setcounter{subsection}{0}

\subsection{Proof of Proposition \ref{lem-for any d}}

\begin{IEEEproof}
 It is clear to see that when $d=1$, $\frac{d}{dx}\left[|1+v'(x)|(1+v'(x))^{-1}\right]=\frac d{dx}[1]=0$. 

By induction, we shall validate \eqref{eqn-nabla=0} for $d\geq2$. Let us denote $V=I+\nabla v^T$, then $\nabla_x^T[|V|V^{-T}]=\nabla_x^T(V^*)^T$, where $V^*$ is the adjugate matrix of $V$ with the element $(V^*)_{ij}=(-1)^{i+j}M_{j,i}$, and $M_{j,i}$ is the $(j,i)$th minor matrix of $V$. Therefore, we have
\begin{align}\label{eqn-nablaV}\notag
	&\nabla_x^T(V^*)^T\\\notag
	=&\begin{pmatrix}
	\frac{\partial}{\partial x_1},\cdots,\frac{\partial}{\partial x_d}
\end{pmatrix}\begin{pmatrix}
	(V^*)_{11}&\cdots&(V^*)_{d1}\\
	\vdots&\ddots&\vdots\\
	(V^*)_{1d}&\cdots&(V^*)_{dd}
\end{pmatrix}\\\notag
	=&\begin{pmatrix}
	\frac{\partial}{\partial x_1},\cdots,\frac{\partial}{\partial x_d}
\end{pmatrix}\begin{pmatrix}
	M_{1,1}&\cdots&(-1)^{1+d}M_{1,d}\\
	\vdots&\ddots&\vdots\\
	(-1)^{d+1}M_{d,1}&\cdots&M_{d,d}
\end{pmatrix}\\
	=&\left[\sum_{i=1}^d(-1)^{i+j}\frac{\partial}{\partial x_i}M_{i,j}\right]_{j=1,\cdots,d}.
\end{align}

Actually, we can show more general statement than \eqref{eqn-nabla=0}. For any $d\geq2$, $i=1,\cdots,d$ and a vector-valued function of $x$ $\bv=(\bv_1,\cdots,\bv_d)\in\mathbb{R}^d$, we have
\begin{align}\label{eqn-Veqn}
	\begin{vmatrix}
		\frac{\partial \bv_1}{\partial x_1}&\cdots&\frac{\partial \bv_{j-1}}{\partial x_1}&\frac{\partial}{\partial x_1}&\frac{\partial \bv_{j+1}}{\partial x_1}&\cdots&\frac{\partial \bv_d}{\partial x_1}\\
		\vdots&\ddots&\vdots&\ddots&\vdots&\ddots&\vdots\\
		\frac{\partial \bv_1}{\partial x_d}&\cdots&\frac{\partial \bv_{j-1}}{\partial x_d}&\frac{\partial}{\partial x_d}&\frac{\partial \bv_{j+1}}{\partial x_d}&\cdots&\frac{\partial \bv_d}{\partial x_d}\end{vmatrix}\equiv0,
\end{align}
where the determinant notation on the left-hand side of \eqref{eqn-Veqn} is defined as the minor expansion along the $j$-th column. 

Equality \eqref{eqn-nabla=0} is just a special case of this statement by letting $\bv_k(x)=x_k+v_k(x)$ in \eqref{eqn-Veqn}. The right-hand side of \eqref{eqn-Veqn} gives exactly the right-hand side of \eqref{eqn-nablaV}.

Now, we shall validate \eqref{eqn-Veqn} for $d\geq2$ by induction.
When $d=2$, direct computation yields that 
\begin{equation}\label{eqn-nabla d=2}
	\begin{vmatrix}
		\frac{\partial}{\partial x_1}&\frac{\partial \bv_2}{\partial x_1}\\
		\frac{\partial}{\partial x_2}&\frac{\partial \bv_2}{\partial x_2}
\end{vmatrix}
=\frac{\partial^2 \bv_2}{\partial x_1\partial x_2}-\frac{\partial^2 \bv_2}{\partial x_2\partial x_1}=0,
\end{equation}
and so does $\begin{vmatrix}
		\frac{\partial \bv_1}{\partial x_1}&\frac{\partial}{\partial x_1}\\
		\frac{\partial \bv_1}{\partial x_2}&\frac{\partial}{\partial x_2}
\end{vmatrix}=0$. 

By induction, suppose \eqref{eqn-Veqn} holds for $d-1$. That is, for any fixed $l,k=1,\cdots,d$, we have
\begin{align}\label{eqn-lemma3}\notag
	0=&\sum_{i\in\{1,\cdots,l-1\}}(-1)^{i+j}\frac{\partial}{\partial x_i}(M^{\bv}_{l,k})_{i,j}\\
	&+\sum_{i\in\{l+1,\cdots,d\}}(-1)^{(i-1)+j}\frac{\partial}{\partial x_i}(M^{\bv}_{l,k})_{i,j},
\end{align}
if $k>j$, and 
\begin{align*}
	0=&\sum_{i\in\{1,\cdots,l-1\}}(-1)^{i+(j-1)}\frac{\partial}{\partial x_i}(M^{\bv}_{l,k})_{i,j}\\
	&+\sum_{i\in\{l+1,\cdots,d\}}(-1)^{(i-1)+(j-1)}\frac{\partial}{\partial x_i}(M^{\bv}_{l,k})_{i,j}=0,
\end{align*}
if $k<j$, where $M^{\bv}_{l,k}$ is the $(l,k)$-th minor matrix of $\nabla v^T$, and $(M^{\bv}_{l,k})_{i,j}$ is the $(i,j)$th minor of $M^{\bv}_{l,k}$.

For any fixed $j=1,\cdots,d$, let us expand $M^{\bv}_{i,j}$ along the $k$th column, $k\in\{1,\cdots,d\}\backslash\{j\}$:
\begin{align}\label{eqn-lemma1}
	&\sum_{i=1}^d(-1)^{i+j}\frac{\partial M_{i,j}}{\partial x_i}\\\notag
=&\left\{\begin{aligned}
	&\sum_{i=1}^d(-1)^{i+j}\frac{\partial}{\partial x_i}
	\left\{\left[\sum_{l\in\{1,\cdots,i-1\}}(-1)^{l+(k-1)}\right.\right.\\
	&\phantom{\sum_{i=1}^d(-1)}\left.\left.+\sum_{l\in\{i+1,\cdots,d\}}(-1)^{(l-1)+(k-1)}\right]\frac{\partial v_k}{\partial x_l}(M_{i,j})_{l,k}\right\},\\
	&\phantom{\sum_{i=1}^d(-1)}\quad\textup{if}\ k>j\\
	&\sum_{i=1}^d(-1)^{i+j}\frac{\partial}{\partial x_i}\left\{\left[\sum_{l\in\{1,\cdots,i-1\}}(-1)^{l+k}\right.\right.\\
	&\phantom{\sum_{i=1}^d(-1)}\left.\left.+\sum_{l\in\{i+1,\cdots,d\}}(-1)^{(l-1)+k}\right]\frac{\partial v_k}{\partial x_l}(M_{i,j})_{l,k}\right\},\\
	&\phantom{\sum_{i=1}^d(-1)}\quad\textup{if}\ k<j
\end{aligned}\right.
\end{align}
Without loss of generality, we shall only examine the right-hand side of \eqref{eqn-lemma1} is identically zero for $k>j$,  since the other case is the same. When $k>j$, we have
\begin{align*}
	&\sum_{i=1}^d(-1)^{i+j}\frac{\partial M^{\bv}_{i,j}}{\partial x_i}\\
\overset{\eqref{eqn-lemma1}}=&\sum_{i=1}^d(-1)^{i+j}\left[\sum_{l\in\{1,\cdots,i-1\}}(-1)^{l+(k-1)}\right.\\
	&\phantom{\sum_{i=1}}\left.+\sum_{l\in\{i+1,\cdots,d\}}(-1)^{(l-1)+(k-1)}\right]\frac{\partial^2\bv_k}{\partial x_i\partial x_l}(M^{\bv}_{i,j})_{l,k}\\	
	&+\sum_{i=1}^d(-1)^{i+j}\left[\sum_{l\in\{1,\cdots,i-1\}}(-1)^{l+(k-1)}\right.\\
	&\phantom{\sum_{i=1}}\left.+\sum_{l\in\{i+1,\cdots,d\}}(-1)^{(l-1)+(k-1)}\right]\frac{\partial \bv_k}{\partial x_l}\frac{\partial}{\partial x_i}(M^{\bv}_{i,j})_{l,k}\\
	&=:I_1+I_2.
\end{align*}
Let us look at the first summation in the bracket in $\uppercase\expandafter{\romannumeral1}_1$:
\begin{align*}
	&\sum_{i=1}^d\sum_{l\in\{1,\cdots,i-1\}}(-1)^{i+j+l+(k-1)}\frac{\partial^2\bv_k}{\partial x_i\partial x_l}(M^{\bv}_{i,j})_{l,k}\\
	=&\sum_{l=1}^d\sum_{i\in\{l+1,\cdots,d\}}(-1)^{i+j+l+(k-1)}\frac{\partial^2\bv_k}{\partial x_i\partial x_l}(M^{\bv}_{i,j})_{l,k}\\
	=&\sum_{i=1}^d\sum_{l\in\{i+1,\cdots,d\}}(-1)^{l+j+i+(k-1)}\frac{\partial^2\bv_k}{\partial x_l\partial x_i}(M^{\bv}_{l,j})_{i,k},
\end{align*}
where the first equality follows by interchanging the order of the two summations. With the fact that $(M^{\bv}_{i,j})_{l,k}=(M^{\bv}_{l,j})_{i,k}$, we have the two summations in the baracket in $I_1$ cancelled out, i.e. $I_1=0$. Meanwhile, by interchanging the order of the summations in $I_2$, one obtain that
\begin{align}\label{eqn-I2}\notag
	I_2
=&\sum_{l=1}^d(-1)^{(l-1)+(k-1)}\frac{\partial \bv_k}{\partial x_l}\\\notag
	&\cdot\left\{\left[\sum_{i\in\{1,\cdots,l-1\}}(-1)^{i+j}+\sum_{i\in\{l+1,\cdots,d\}}(-1)^{(i-1)+j}\right]\right.\\
	&\phantom{aaa}\left.\cdot\frac{\partial}{\partial x_i}(M^{\bv}_{l,k})_{i,j}\right\}\overset{\eqref{eqn-lemma3}}=0,
\end{align}
since the terms in the brace equal zero, by induction hypothesis. The case $k<j$ can be argued similarly to verify that \eqref{eqn-Veqn} holds.
\end{IEEEproof}

\subsection{Proof of Proposition \ref{prop-2.2}}
\begin{IEEEproof}
We move the term with $h$ in \eqref{eqn-K} to the left-hand side and divide by $p^2$ throughout:
\begin{align}\label{eqn-nablah}\notag
	-\nabla^Th=&\frac1pK^T\nabla^2 p+\frac1p\nabla^Tp(\nabla K^T)^T\\\notag
	&+\left[\textup{tr}\left(\frac{\partial(\nabla K^T)}{\partial x_i}\right)\right]_{i=1,\cdots,d}-\frac1{p^2}(\nabla^Tp)K\nabla^Tp\\
=:I&II_1+III_2+III_3+III_4.
\end{align}
Notice that 
\begin{align}\label{eqn-nabla2lnp}\notag
	\nabla^2(\log p)=&\nabla\left(\nabla^T(\log p)\right)=\nabla\left(\frac1p\nabla^Tp\right)\\\notag
=&-\frac1{p^2}\nabla p\nabla^Tp+\frac1p\nabla^2p\\
	=&-\nabla(\log p)\nabla^T(\log p)+\frac1p\nabla^2p,
\end{align}
then 
\begin{align}\label{eqn-III1}
	III_1=&K^T\left[\nabla^2(\log p)+\nabla(\log p)\nabla^T(\log p)\right].
\end{align}
With the fact that for any $i=1,\cdots,d$,
\begin{align}\label{eqn-traceK''}\notag
	\textup{tr}\left(\frac{\partial (\nabla K^T)}{\partial x_i}\right)
=&\textup{tr}\left[\frac{\partial}{\partial x_i}\begin{pmatrix}
	\frac{\partial K_1}{\partial x_1}&\cdots&\frac{\partial K_d}{\partial x_1}\\
	\vdots&\ddots&\vdots\\
	\frac{\partial K_1}{\partial x_d}&\cdots&\frac{\partial K_d}{\partial x_d}
\end{pmatrix}\right]\\
	=&\frac{\partial }{\partial x_i}(\nabla^TK),
\end{align}
we have
\begin{equation}\label{eqn-III3}
	III_3=\nabla^T(\nabla^T K).
\end{equation}
Lastly, it is easy to see that
\begin{align}\label{eqn-III2}
	III_2=&\nabla^T(\log p)(\nabla K^T)^T,\\\label{eqn-III4}
	III_4=&-\nabla^T(\log p)K\nabla^T(\log p).
\end{align}
Substituting \eqref{eqn-III1}, \eqref{eqn-III3}-\eqref{eqn-III4} back to \eqref{eqn-nablah}, we obtain that
\begin{equation}\label{eqn-nablah1}
	-\nabla^Th=K^T\nabla^2(\log p)+\nabla^T(\log p)(\nabla K^T)^T+\nabla^T(\nabla^T K),
\end{equation}
since the second term of $III_1$ cancels out with $III_4$. Let us verify that the right-hand side of \eqref{eqn-nablah1} is $\nabla^T\left(\frac1p\nabla^T(pK)\right)$ by direct computation:
\begin{align}\label{eqn-K-1}
	&\nabla^T\left(\frac1p\nabla^T(pK)\right)\\\notag
=&\nabla^T\left(\frac{\nabla^TpK}p+\nabla^TK\right)
=\nabla^T\left(\nabla^T(\log p)K+\nabla^TK\right)\\\notag
=&K^T[\nabla(\nabla^T(\log p))]^T+\nabla^T(\log p)(\nabla K^T)^T+\nabla^T(\nabla^TK),
\end{align}
where $\nabla(\nabla^T(\log p))^T=\nabla^2(\log p)$ due to its symmetry and the third equality follows from the fact that, for any two column vectors $a$, $b\in\R^d$, 
\begin{equation}\label{eqn-nablaab}
	\nabla^T(a^Tb)=b^T(\nabla a^T)^T+a^T(\nabla b^T)^T.
\end{equation}

To derive \eqref{eqn-u-2}, we look at the terms in \eqref{eqn-u}. It contains the same terms as those $III_1-III_4$ in \eqref{eqn-nablah} with $K$ replacing by $u$. Therefore, by moving all these terms with $u$ to the left-hand side and dividing by $p^2$ throughout, we have
\begin{align}\label{eqn-u-1}\notag
	&-\nabla^T\left(\frac1p\nabla^T(pu)\right)\\\notag
\overset{\eqref{eqn-u},\eqref{eqn-K-1}}=&-h\nabla^Th+K^T\nabla^2h\\\notag
	&+\frac1{2p}K^T\left(\nabla^T\otimes\left(\nabla^2p\right)\right)(I_d\otimes K)\\\notag
	&+\frac1pK^T\nabla^2p(\nabla K^T)^T+\frac1p(\nabla^Tp)K\nabla^Th\\\notag
	&+\nabla^Th(\nabla K^T)^T+\frac1p(\nabla^Tp)K\left[\textup{tr}\left(\frac{\partial(\nabla K^T)}{\partial x_i}\right)\right]_{i=1,\cdots,d}\\\notag
	&-\left[\textup{tr}\left(\nabla K^T\frac{\partial(\nabla K^T)}{\partial x_i}\right)\right]_{i=1,\cdots,d}\\\notag
	&-\frac1{2p^2}K^T(\nabla^2p)K\nabla^Tp\\
	=:&-h\nabla^Th+IV_1+IV_2+\cdots+IV_8.
\end{align} 
In the sequel, we shall simplify the terms $IV_1-IV_8$ one-by-one into the derivatives of $\log p$ and $K$:
\begin{align*}
	IV_1=&K^T\nabla(\nabla^Th)\\
\overset{\eqref{eqn-nablah1}}=&-K^T\nabla\left[K^T\nabla^2(\log p)+\nabla^T(\log p)(\nabla K^T)^T\right.\\
	&\phantom{-K^T\nabla[]}\left.+\nabla^T(\nabla^TK)\right]\\
	=&-K^T(\nabla K^T)\nabla^2(\log p)\\
	&-K^T(I_d\otimes K^T)(\nabla\otimes\nabla^2(\log p))\\
	&-K^T\nabla^2(\log p)(\nabla K^T)^T\\
	&-K^T(I_d\otimes\nabla^T(\log p))(\nabla\otimes(\nabla K^T)^T)\\
	&-K^T\nabla^2(\nabla^TK)\\
	=:&IV_{1,1}+IV_{1,2}+IV_{1,3}+IV_{1,4}+IV_{1,5},
\end{align*}
where the second equality follows from the fact that for a column vector $a\in\R^d$ and a $d\times d$ matrix $A$, 
\begin{align}\label{eqn-nablaTaA}
	\nabla(a^TA)=(\nabla a^T)A+(I_d\otimes a^T)(\nabla\otimes A).
\end{align}
\begin{align*}
	IV_2=&\frac12\left[K^T\frac1p\frac{\partial}{\partial x_i}(\nabla^2 p)K\right]_{i=1,\cdots,d}\\
\overset{\eqref{eqn-nabla2lnp}}=&\frac12\left[K^T\frac1p\frac{\partial}{\partial x_i}(p\nabla^2(\log p)\right.\\
	&\phantom{\frac12aa}\left.+p\nabla(\log p)\nabla^T(\log p))K\right]_{i=1,\cdots,d}\\
	=&\frac12\left\{K^T\left[\frac{\partial\log p}{\partial x_i}(\nabla^2(\log p)+\nabla(\log p)\nabla^T(\log p))\right.\right.\\
	&\phantom{\frac12K^T[][]}+\frac{\partial}{\partial x_i}\nabla^2(\log p)+\frac{\partial}{\partial x_i}(\nabla(\log p))\nabla^T(\log p)\\
	&\phantom{\frac12K^T[][]}\left.\left.+\nabla(\log p)\frac{\partial}{\partial x_i}(\nabla^T(\log p))\right]K\right\}_{i=1,\cdots,d}\\
	=&\frac12K^T\left[\nabla^T(\log p)\otimes\nabla^2(\log p)\right](I_d\otimes K)\\
	&+\frac12K^T\left[\nabla^T(\log p)\otimes\nabla(\log p)\nabla^T(\log p)\right](I_d\otimes K)\\
	&+\frac12K^T\left[\nabla^T\otimes\nabla^2(\log p)\right](I_d\otimes K)\\
	&+\frac12K^T\nabla^2(\log p)\nabla^T(\log p)K\\
	&+\frac12K^T\nabla(\log p)(\nabla^T\otimes\nabla^T(\log p))(I_d\otimes K)\\
=:&IV_{2,1}+IV_{2,2}+IV_{2,3}+IV_{2,4}+IV_{2,5}.
\end{align*}
\begin{align*}
	IV_3\overset{\eqref{eqn-nabla2lnp}}=&K^T\left[\nabla^2(\log p)+\nabla(\log p)\nabla^T(\log p)\right](\nabla K^T)^T\\
	=&K^T\nabla^2(\log p)(\nabla K^T)^T\\
	&+K^T\nabla(\log p)\nabla^T(\log p)(\nabla K^T)^T\\
	=:&IV_{3,1}+IV_{3,2}.
\end{align*}
\begin{align*}
	IV_4\overset{\eqref{eqn-nablah1}}=&-\nabla^T(\log p)K\\
	&\phantom{aa}\cdot\left[K^T\nabla^2(\log p)+\nabla^T(\log p)(\nabla K^T)^T\right.\\
	&\phantom{aaaa}\left.+\nabla^T(\nabla^TK)\right]\\
	=&-\nabla^T(\log p)KK^T\nabla^2(\log p)\\
	&-\nabla^T(\log p)K\nabla^T(\log p)(\nabla K^T)^T\\
	&-\nabla^T(\log p)K\nabla^T(\nabla^TK)\\
	=:&IV_{4,1}+IV_{4,2}+IV_{4,3}.
\end{align*}
\begin{align*}
	IV_5\overset{\eqref{eqn-nablah1}}
	=&-\left[K^T\nabla^2(\log p)+\nabla^T(\log p)(\nabla K^T)^T\right.\\
	&\phantom{-aa}+\nabla^T(\nabla^TK)]\cdot(\nabla K^T)^T\\
	=&-K^T\nabla^2(\log p)(\nabla K^T)^T-\nabla^T(\log p)\left[(\nabla K^T)^T\right]^2\\
	&-\nabla^T(\nabla^TK)(\nabla K^T)^T\\
	=:&IV_{5,1}+IV_{5,2}+IV_{5,3}.
\end{align*}
\begin{align*}
	IV_6\overset{\eqref{eqn-traceK''}}=&\nabla^T(\log p)K\nabla^T(\nabla^TK).
\end{align*}
\begin{align*}
	IV_7=\frac12\nabla^T\left\{\textup{tr}\left[(\nabla K^T)(\nabla K^T)^T\right]\right\}.
\end{align*}
In fact, for fixed $l=1,\cdots,d$:
\begin{align}\label{eqn-IV7}\notag
	&\textup{tr}\left[\nabla K^T\frac{\partial(\nabla K^T)}{\partial x_l}\right]\\\notag
	=&\sum_{i,j=1}^d\frac{\partial K_j}{\partial x_i}\frac{\partial^2 K_i}{\partial x_l\partial x_j}\\\notag
	=&\frac{\partial}{\partial x_l}\sum_{i,j=1}^d\left[\frac{\partial K_j}{\partial x_i}\frac{\partial K_i}{\partial x_j}\right]-\sum_{i,j=1}^d\left[\frac{\partial}{\partial x_l}\left(\frac{\partial K_j}{\partial x_i}\right)\frac{\partial K_i}{\partial x_j}\right]\\
	=&\frac{\partial}{\partial x_l}\textup{tr}\left[(\nabla K^T)(\nabla K^T)^T\right]-\textup{tr}\left[\nabla K^T\frac{\partial}{\partial x_l}(\nabla K^T)\right].
\end{align}
Consequently, $IV_7$ follows by combining the similar terms in \eqref{eqn-IV7} and dividing by $2$.
\begin{align*}
	IV_8=&-\frac12K^T\left(\frac1p\nabla^2p\right)K\nabla^T(\log p)\\
\overset{\eqref{eqn-nabla2lnp}}=&-\frac12K^T\left[\nabla^2(\log p)+\nabla(\log p)\nabla^T(\log p)\right]K\nabla^T(\log p)\\
	=&-\frac12K^T\nabla^2(\log p)K\nabla^T(\log p)\\
	&-\frac12K^T\nabla(\log p)\nabla^T(\log p)K\nabla^2(\log p)\\
	=:&IV_{8,1}+IV_{8,2}.
\end{align*}

After close examination to the terms in $IV_1-IV_8$, we obtain the following identities:
\begin{enumerate}
	\item $IV_{1,3}=-IV_{3,1}$;
	\item $IV_{3,2}=-IV_{4,2}$, due to the fact that $K^T\nabla(\log p)=\nabla^T(\log p)K$;
	\item $IV_{4,3}=-IV_6$;
	\item $IV_{2,1}=-IV_{8,1}$ and $IV_{2,2}=-IV_{8,2}$;
	\item $IV_{2,4}+IV_{2,5}+IV_{4,1}=0$;
	\item $IV_{1,2}+IV_{2,3}+IV_{5,1}=-\frac12\nabla^T\left[\frac1p(K^T\nabla^2pK)\right]+\frac12\nabla^T\left[K^T\nabla(\log p)\nabla^T(\log p)K\right]$;
	\item $IV_{1,1}+IV_{1,4}+IV_{5,2}=-\nabla^T\left[K^T(\nabla K^T)\nabla(\log p)\right]$;
	\item $IV_{1,5}+IV_{5,3}=-\nabla^T\left[\nabla^T(\nabla^TK)K\right]$.
\end{enumerate}

The necessary proofs of these identities are postponed to the end of the proof. Equation \eqref{eqn-u-2} follows by substituting these identity 1)-8) back to $IV_1-IV_8$. 
\end{IEEEproof}

In the rest of this subsection, we shall detail the proof of identity 4)-8).

\begin{IEEEproof}[Proof of 4)]
, which holds due to the following identities
\begin{align*}
	&K^T\left[\nabla^T(\log p)\otimes\nabla^2(\log p)\right](I_d\otimes K)\\
=&(1\otimes K^T)\left\{\nabla^T(\log p)\otimes[\nabla^2(\log p)K]\right\}\\
=&\nabla^T(\log p)\otimes[K^T\nabla^(\log p)K]=K^T\nabla^2(\log p)K\nabla^T(\log p),
\end{align*}
and 
\begin{align*}
	&K^T\left\{\nabla^T(\log p)\otimes[\nabla(\log p)\nabla^T(\log p)]\right\}(I_d\otimes K)\\
	=&\left(1\otimes K^T\right)\left\{\nabla^T(\log p)\otimes[\nabla(\log p)\nabla^T(\log p)K]\right\}\\
	=&\nabla^T(\log p)\otimes\left[K^T\nabla(\log p)\nabla^T(\log p)K\right]\\
	=&K^T\nabla(\log p)\nabla^T(\log p)K\nabla^T(\log p),
\end{align*}
where the last equalities in both identities are obtained by noting that $K^T\nabla^2(\log p)K$ and $K^T\nabla(\log p)\nabla^T(\log p)K$ are scalars.
\end{IEEEproof}

\begin{IEEEproof}[Proof of 5)]
Notice that $K^T\nabla(\log p)=\nabla^T(\log p)K$ is a scalar. Divide by this scalar through out these three terms, it yields that
\begin{align*}
	&\frac1{K^T\nabla(\log p)}\left(IV_{2,4}+IV_{2,5}+IV_{4,1}\right)\\
=&\frac12K^T\nabla^2(\log p)+\frac12(\nabla^T\otimes\nabla^T(\log p))(I_d\otimes K)\\
	&-K^T\nabla^2(\log p)\\
=&\frac12\left[(\nabla^T\otimes\nabla^T(\log p))(I_d\otimes K)-K^T\nabla^2(\log p)\right]
	=0,
\end{align*}
since 
\begin{align*}
	&\left[(\nabla^T\otimes\nabla^T(\log p))(I_d\otimes K)\right]_j
=\left[\frac{\partial\nabla^T(\log p)}{\partial x_j}K\right]_{j}\\
	=&\left[\sum_{i=1}^d\frac{\partial^2(\log p)}{\partial x_j\partial x_i}K_i\right]_j
	=\left[K^T\nabla^2(\log p)\right]_j,
\end{align*}
for $j=1,\cdots,d$.
\end{IEEEproof}

\begin{IEEEproof}[Proof of 6)]
	Let us begin with the direct computation $\nabla^T\left[\frac1p(K^T\nabla^2pK)\right]$:
\begin{align}\label{eqn-nabla1pknabla2pk}\notag
	&\nabla^T\left[\frac1p(K^T\nabla^2pK)\right]\\\notag
	\overset{\eqref{eqn-nabla2lnp}}=&\nabla^T\left[K^T\nabla(\log p)\nabla^T(\log p)K\right]+\nabla^T\left[K^T\nabla^2(\log p)K\right]\\\notag
	=&\nabla^T\left[K^T\nabla(\log p)\nabla^T(\log p)K\right]+\left[\nabla^2(\log p)K\right]^T(\nabla K^T)^T\\\notag
	&+K^T\left\{\nabla\left[\nabla^2(\log p)K\right]^T\right\}^T\\\notag
	=&\nabla^T\left[K^T\nabla(\log p)\nabla^T(\log p)K\right]+\left[\nabla^2(\log p)K\right]^T(\nabla K^T)^T\\\notag
	&+K^T\left\{\nabla\left[K^T\nabla^2(\log p)\right]\right\}^T\\\notag
	\overset{\eqref{eqn-nablaTaA}}=&\nabla^T\left[K^T\nabla(\log p)\nabla^T(\log p)K\right]+\left[\nabla^2(\log p)K\right]^T(\nabla K^T)^T\\\notag
	&+K^T\left[\nabla K^T\nabla^2(\log p)+(I_d\otimes K^T)(\nabla\otimes\nabla^2(\log p))\right]^T\\\notag
	=&\nabla^T\left[K^T\nabla(\log p)\nabla^T(\log p)K\right]\\\notag
			&+2K^T\nabla^2(\log p)(\nabla K^T)^T\\
	&+K^T\left[\nabla^T\otimes \nabla^2(\log p)\right](I_d\otimes K).
\end{align}
Moreover, we notice that the terms $IV_{1,2}$ and $IV_{2,3}$ can be combined, since
\begin{align}\label{eqn-IV12 and IV23}
	(I_d\otimes K^T)(\nabla\otimes \nabla^2(\log p))=(\nabla^T\otimes \nabla^2(\log p))(I_d\otimes K).
\end{align}
In fact, the left-hand side of \eqref{eqn-IV12 and IV23} is
\begin{align*}
	(I_d\otimes K^T)(\nabla\otimes \nabla^2(\log p))
	=&\begin{bmatrix}
	K^T\frac{\partial \nabla^2(\log p)}{\partial x_1}\\
	\vdots\\
	K^T\frac{\partial \nabla^2(\log p)}{\partial x_d}
\end{bmatrix},
\end{align*}
where the $(l,j)$th entry of 
\begin{equation}\label{eqn-p'''-l}
	\left[(I_d\otimes K^T)(\nabla\otimes \nabla^2(\log p))\right]_{lj}=\sum_{i=1}^dK_i\frac{\partial^3(\log p)}{\partial x_l\partial x_i\partial x_j},
\end{equation}
while the right-hand side of \eqref{eqn-IV12 and IV23} is
\begin{align*}
	(\nabla^T\otimes \nabla^2(\log p))(I_d\otimes K)=\begin{bmatrix}
	\frac{\partial \nabla^2(\log p)}{\partial x_1}K,\cdots,\frac{\partial \nabla^2(\log p)}{\partial x_d}K
\end{bmatrix},
\end{align*}
where the $(i,l)$th entry of
\begin{equation}\label{eqn-p'''-r}
	\left[(\nabla^T\otimes \nabla^2(\log p))(I_d\otimes K)\right]_{il}=\sum_{j=1}^d\frac{\partial^3(\log p)}{\partial x_l\partial x_i\partial x_j}K_j.
\end{equation}
Equality \eqref{eqn-IV12 and IV23}  immediately follows by replacing $(i,j,l)$ in \eqref{eqn-p'''-r} into $(l,i,j)$ in \eqref{eqn-p'''-l}. 

Combining \eqref{eqn-nabla1pknabla2pk} and \eqref{eqn-IV12 and IV23}, we finished the proof of (6). 
\end{IEEEproof}
\begin{IEEEproof}[Proof of 7)]
	Through direct computations of $\nabla^T\left[K^T\nabla K^T\nabla(\log p)\right]$, we have
	\begin{align}\label{eqn-nablaknablaktnablalnp}\notag
		&\nabla^T\left[K^T\nabla K^T\nabla(\log p)\right]\\\notag
		\overset{\eqref{eqn-nablaab}}=&\nabla^T(\log p)\left[\nabla(K^T\nabla K^T)\right]^T+K^T\nabla K^T\left[\nabla^2(\log p)\right]^T\\\notag
		\overset{\eqref{eqn-nablaTaA}}=&\nabla^T(\log p)\left\{(\nabla K^T)^2+(I_d\otimes K^T)\left[\nabla\otimes(\nabla K^T)\right]\right\}^T\\\notag
	&+K^T\nabla K^T\left[\nabla^2(\log p)\right]^T\\\notag
		=&\nabla^T(\log p)\left[(\nabla K^T)^T\right]^2\\\notag
		&+\nabla^T(\log p)\left[\nabla^T\otimes(\nabla K^T)^T\right](I_d\otimes K)\\
	&+K^T\nabla K^T\nabla^2(\log p).
	\end{align}
It is easy to check that 
\begin{align}\label{eqn-IV14}\notag
	-IV_{1,4}=&K^T(I_d\otimes\nabla^T(\log p))\left[\nabla\otimes(\nabla K^T)^T\right]\\
	=&\nabla^T(\log p)\left[\nabla^T\otimes(\nabla K^T)^T\right](I_d\otimes K),
\end{align}
since the $j$th component of the middle term of \eqref{eqn-IV14} is
\begin{align*}
	&\left\{K^T(I_d\otimes\nabla^T(\log p))\left[\nabla\otimes(\nabla K^T)^T\right]\right\}_j\\
	=&\sum_{i=1}^dK_i\left\{\left[I_d\otimes\nabla^T(\log p)\right]\left[\nabla\otimes(\nabla K^T)^T\right]\right\}_{ij}\\
	=&\sum_{i=1}^dK_i\sum_{l=1}^d\frac{\partial \log p}{\partial x_l}\frac{\partial \left[(\nabla K^T)^T\right]_{lj}}{\partial x_i}\\
	=&\sum_{i=1}^dK_i\sum_{l=1}^d\frac{\partial \log p}{\partial x_l}\frac{\partial^2K_l}{\partial x_i\partial x_j},
\end{align*}
while the $k$th component of the right-hand side is 
\begin{align*}
	&\left\{\nabla^T(\log p)\left[\nabla^T\otimes(\nabla K^T)^T\right](I_d\otimes K)\right\}_k\\
	=&\sum_{l=1}^d\frac{\partial\log p}{\partial x_l}\left\{\left[\nabla^T\otimes(\nabla K^T)^T\right](I_d\otimes K)\right\}_{lk}\\
	=&\sum_{l=1}^d\frac{\partial\log p}{\partial x_l}\sum_{j=1}^d\frac{\partial \left[(\nabla K^T)^T\right]_{lj}}{\partial x_k}K_j\\
	&=\sum_{l=1}^d\frac{\partial\log p}{\partial x_l}\sum_{j=1}^d\frac{\partial^2K_l}{\partial x_j\partial x_k}K_j.
\end{align*}
Item (7) follows immediately from \eqref{eqn-nablaknablaktnablalnp} and \eqref{eqn-IV14}.
\end{IEEEproof}
\begin{IEEEproof}[Proof of 8)]
	Item (8) follows immediately by directly computing $\nabla^T\left[\nabla^T(\nabla^TK)K\right]$:
	\begin{align*}
	\nabla^T\left[\nabla^T(\nabla^TK)K\right]
=&\nabla^T\left\{\left[\nabla^T(\nabla^TK)\right]K\right\}\\
\overset{\eqref{eqn-nablaab}}=&K^T\nabla^2(\nabla^TK)+\nabla^T(\nabla^TK)(\nabla K^T)^T.
\end{align*}
\end{IEEEproof}

\subsection{Sobolev embedding theorem}
\renewcommand{\thetheorem}{C.\arabic{theorem}}

Let $W^{k,p}(\R^d;w)$ denote the Sobolev space consisting of all real-valued functions on $\R^d$ whose first $k$ weak derivatives are functions in $L^p(\R^d;w)$. Here $k$ is a non-negative integer and $1\leq p<\infty$. 
\begin{theorem}[\cite{E,GT}]\label{thm-Sobolev}
	 If $k>l$ and $1\leq p<q<\infty$ are two real numbers such that $(k-l)p<n$ and $\frac1p-\frac kd=\frac1q-\frac ld$, then $W^{k,p}(\R^d;w)\subset W^{l,q}(\R^d;w)$.
\end{theorem} 

\subsection{Proof of Lemma \ref{lem-1}}

\begin{IEEEproof}
	We proceed by induction. For the case $k=1$, it can be easily verified by differentiating \eqref{eqn-K-2} with respect to $x_{\sigma_1(1)}$, for any $\sigma_1\in\{1,\cdots,d\}$:
\begin{align*}
	\textup{LHS}=&-\nabla^T\left[p\nabla\left(\frac{\partial \phi}{\partial x_{\sigma_1(1)}}\right)\right]-\nabla^T\left[\frac{\partial p}{\partial x_{\sigma_1(1)}}\nabla\phi\right]\\
\overset{\eqref{eqn-K-2}}=&-\nabla^T\left[p\nabla\left(\frac{\partial \phi}{\partial x_{\sigma_1(1)}}\right)\right]-\nabla^T\left[\frac{\partial(\log p)}{\partial x_{\sigma_1(1)}}\right]p\nabla\phi\\
	&+\frac{\partial(\log p)}{\partial x_{\sigma_1(1)}}(h-\hat{h})p,
\end{align*}
while 
\begin{align*}
	\textup{RHS}=&\frac {\partial h}{\partial x_{\sigma_1(1)}}p+(h-\hat{h})\frac{\partial(\log p)}{\partial x_{\sigma_1(1)}}p.
\end{align*}
The last terms on the both sides are cancelled out with each other. Thus, we have
\begin{align*}
	-\nabla^T\left[p\nabla\left(\frac{\partial \phi}{\partial x_{\sigma_1(1)}}\right)\right]
=&\nabla^T\left[\frac{\partial(\log p)}{\partial x_{\sigma_1(1)}}\right]p\nabla\phi+\frac {\partial h}{\partial x_{\sigma_1(1)}}p\\
	=:&G^1_{\sigma_1}p,
\end{align*}
with $G^1_{\sigma_1}$ defined in \eqref{eqn-G1}. Suppose \eqref{eqn-DkeqnK} holds for $k=l-1$, i.e.
\begin{align}\label{eqn-Dl-1}
	-\nabla^T\left[p\nabla\left(\nabla_{\sigma_{l-1}}^{l-1}\phi\right)\right]
	=G_{\sigma_{l-1}}^{l-1}p,
\end{align}
where $\sigma_{l-1}\in\{1,\cdots,d\}^{l-1}$, then we shall validate it for $k=l$. Differentiating \eqref{eqn-Dl-1} with respect to $x_{\sigma(l)}$, $\sigma(l)\in\{1,\cdots,d\}$ yields that
\begin{align*}
	\textup{LHS}\overset{\eqref{eqn-Dl-1}}
	=&-\nabla^T\left[p\nabla\left(\nabla_{(\sigma_{l-1},\sigma(l))}^{l}\phi\right)\right]
\\
	&-\nabla^T\left[\frac{\partial(\log p)}{\partial x_{\sigma(l)}}\right]p\nabla\left(\nabla_{\sigma_{l-1}}^{l-1}\phi\right)+\frac{\partial(\log p)}{\partial x_{\sigma(l)}}G_{\sigma_{l-1}}^{l-1}p,
\end{align*}
while
\begin{align*}
	\textup{RHS}=&\frac{\partial G^{l-1}_{\sigma_{l-1}}}{\partial x_{\sigma(l)}}p+ G^{l-1}_{\sigma_{l-1}}\frac{\partial p}{\partial x_{\sigma(l)}}.
\end{align*}
Similarly, the last terms on both sides are cancelled out. Thus, we have
\begin{align*}
	&-\nabla^T\left[p\nabla\left(\nabla_{(\sigma_{l-1},\sigma(l))}^{l}\phi\right)\right]\\
=&\nabla^T\left[\frac{\partial(\log p)}{\partial x_{\sigma(l)}}\right]p\nabla\left[\nabla^{l-1}_{\sigma_{l-1}}\phi\right]+\frac{\partial G^{l-1}_{\sigma_{l-1}}}{\partial x_{\sigma(l)}}p,
\end{align*}
which is exactly \eqref{eqn-DkeqnK} with $G_{\sigma_k}^k$ defined in \eqref{eqn-Gk} with $\sigma_k=(\sigma_{l-1},\sigma(l)$. Next, we compute $\frac{\partial G^{l-1}_{\sigma_{l-1}}}{\partial x_{\sigma(l)}}$, for $l\geq2$:
\begin{align*}
	&\frac{\partial G^{l-1}_{\sigma_{l-1}}}{\partial x_{\sigma(l)}}\\
\overset{\eqref{eqn-Gk}}=&\frac\partial{\partial x_{\sigma(l)}}\left[\nabla^T\left(\frac{\partial(\log p)}{\partial x_{\sigma_{l-1}(l-1)}}\right)\nabla\left(\nabla^{l-2}_{\sigma_{l-1}|_{1:l-2}}\phi\right)\right]\\
	&+\frac{\partial^2G^{l-2}_{\sigma_{l-1}|_{1:l-2}}}{\partial x_{\sigma(l)}\partial x_{\sigma_{l-1}(l-1)}}\\
	=&\nabla_{\sigma(l)}\left[\nabla^T\left(\frac{\partial(\log p)}{\partial x_{\sigma_{l-1}(l-1)}}\right)\nabla\left(\nabla^{l-2}_{\sigma_{l-1}|_{1:l-2}}\phi\right)\right]\\
	&+\sum_{k=3}^{l-1}\nabla^{l-k+1}_{\sigma(l),\sigma_{l-1}|_{k:l-1}}\left[\nabla^T\left(\frac{\partial(\log p)}{\partial x_{\sigma_{l-1}(k-1)}}\right)\right.\\
	&\phantom{+\sum_{k=2}^{l-1}\nabla^k_{\sigma(l),\sigma_{l-1}|_{k:l-1}}aa}\left.\cdot\nabla\left(\nabla^{k-2}_{\sigma_{l-1}|_{1:k-2}}\phi\right)\right]\\
	&+\nabla^{l-1}_{(\sigma(l),\sigma_{l-1}|_{2:l-1})}G^{1}_{\sigma_{l-1}(1)},
\end{align*}
with 
\begin{align*}
	&\nabla^{l-1}_{(\sigma(l),\sigma_{l-1}|_{2:l-1})}G^{1}_{\sigma_{l-1}(1)}\\
\overset{\eqref{eqn-G1}}=&\nabla^{l-1}_{(\sigma(l),\sigma_{l-1}|_{2:l-1})}\left[\nabla^T\left(\frac{\partial(\log p)}{\partial x_{\sigma_{l-1}(1)}}\right)\nabla\phi\right]+\nabla^l_{\sigma(l),\sigma_{l-1}}h.
\end{align*}
Equation \eqref{eqn-parialGk} follows immediately by denoting $\sigma_l=(\sigma_{l-1},\sigma(l))\in\{1,\cdots,d\}^l$.
\end{IEEEproof}

\end{document}